\documentclass[11pt]{article}
\usepackage[english]{babel}
\usepackage{indentfirst}
\usepackage{latexsym}
\usepackage[usenames,dvipsnames]{color}
\usepackage[colorlinks=true, bookmarksnumbered=true, bookmarksopen=true,
bookmarksopenlevel=3, pdfstartview=FitH, linkcolor=cyan, pdfmenubar=true,
pdftoolbar=true, bookmarks=true,citecolor=cyan, urlcolor=magenta,
filecolor=cyan,menucolor=black,plainpages=false,pdfpagelabels]{hyperref}
\usepackage[lmargin=2cm,tmargin=2cm,bmargin=2cm,rmargin=2cm]{geometry}
\usepackage{amsbsy,amsmath,amsfonts,amssymb,amsthm}
\usepackage{times}
\usepackage{graphicx}
\usepackage[utf8]{inputenc}
\usepackage[T1]{fontenc}

\numberwithin{equation}{section}
\newtheorem{prop}{Proposition}[section]

\newtheorem{theorem}[prop]{Theorem}
\newtheorem{cor}[prop]{Corollary}
\newtheorem{remark}[prop]{Remark}
\newtheorem{lemma}[prop]{Lemma}

\newcommand\reallywidehat[1]{\savestack{\tmpbox}{\stretchto{  \scaleto{\scalerel*[\widthof{\ensuremath{#1}}]{\kern-.6pt\bigwedge\kern-.6pt}{\rule[-\textheight/2]{1ex}{\textheight}}  }{\textheight}}{0.5ex}}\stackon[1pt]{#1}{\tmpbox}}
\def\fin { \vskip 0pt \hfill $\diamond$ \vskip 12pt}

\begin{document}

\title{Global well-posedness and asymptotic behavior for
Navier-Stokes-Coriolis equations in homogeneous Besov spaces}
\author{{Lucas C. F. Ferreira$^{1}$}{\thanks{{Corresponding author. }\newline
{E-mail addresses: lcff@ime.unicamp.br (L.C.F. Ferreira),
vladimirangulo01@gmail.com (V. Angulo-Castillo).}\newline
{LCF Ferreira was supported by FAPESP and CNPq, Brazil. V. Angulo-Castillo
was supported by CNPq, Brazil.}}}, \ Vladimir Angulo-Castillo$^{2}$ \\
{\small $^{1,2}$ Universidade Estadual de Campinas, IMECC-Departamento de
Matem\'{a}tica,}\\
{\small Rua S\'{e}rgio Buarque de Holanda, CEP 13083-859, Campinas, SP,
Brazil.}}
\date{}
\maketitle

\begin{abstract}
We are concerned with the $3$D-Navier-Stokes equations with Coriolis force.
Existence and uniqueness of global solutions in homogeneous Besov spaces are
obtained for large speed of rotation. In the critical case of the
regularity, we consider a suitable initial data class whose definition is
based on the Stokes-Coriolis semigroup and Besov spaces. Moreover, we
analyze the asymptotic behavior of solutions in that setting as the speed of
rotation goes to infinity. \medskip


{\small \medskip\noindent\textbf{Key:} Navier-Stokes equations; Coriolis
force; Global well-posedness; Asymptotic behavior; Besov spaces}
\end{abstract}

\renewcommand{\abstractname}{Abstract}

\section{Introduction}

In this paper we are concerned with the incompressible Navier-Stokes
equations in the rotational framework
\begin{equation}
\left\{ \begin{aligned} & \frac{\partial u}{\partial t}-\Delta u+\Omega
e_{3}\times u+ \left(u\cdot \nabla\right) u + \nabla p = 0 \text{ in }\
\mathbb{R}^{3}\times (0,\infty)\\ & \nabla\cdot u = 0 \text{ in }\
\mathbb{R}^{3}\times (0,\infty)\\ & u(x,0) = u_{0}(x) \text{ in }\
\mathbb{R}^{3} \end{aligned}\right. ,  \label{NSC}
\end{equation}%
where $u=u(x,t)=\left( u_{1}(x,t),u_{2}(x,t),u_{3}(x,t)\right) $ and $%
p=p(x,t)$ stand for the velocity field and the pressure of the fluid,
respectively. The initial data $u_{0}=\left(
u_{0,1}(x),u_{0,2}(x),u_{0,3}(x)\right)$ satisfies the divergence-free
condition $\nabla \cdot u_{0}=0$. The letter $\Omega \in \mathbb{R}$
represents the Coriolis parameter while its modulus $|\Omega |$ is the speed
of rotation around the vertical vector $e_{3}=(0,0,1)$. For more details
about the physical model, we refer the reader to the book \cite%
{Chemin-Gallagher2006}. Here, we will use the same notation for spaces of
scalar and vector functions, e.g., we write $u_{0}\in H^{s}$ instead of $%
u_{0}\in (H^{s})^{3}$.

Invoking Duhamel's principle, the system (\ref{NSC}) can be converted to the
integral equation (see e.g. \cite{HieberShibata2010})
\begin{equation}
u(t)=T_{\Omega }(t)u_{0}-\mathfrak{B}(u,u)(t),  \label{integralform}
\end{equation}%
where the bilinear operator $\mathfrak{B}$ is defined by
\begin{equation}
\mathfrak{B}(u,v)(t)=\int_{0}^{t}T_{\Omega }(t-\tau )\mathbb{P}\nabla \cdot
(u\otimes v)(\tau )\ d\tau .  \label{bili-aux-1}
\end{equation}%
In (\ref{bili-aux-1}), $\mathbb{P}=(\delta _{i,j}+R_{i}R_{j})_{1\leq i,j\leq
3}$ is the Leray-Helmholtz projector, $\{R_{i}\}_{1\leq i\leq 3}$ are the
Riesz transforms, and $T_{\Omega }(\cdot )$ stands for the semigroup
corresponding to the linear part of (\ref{NSC}) (Stokes-Coriolis semigroup).
More explicitly, we have that%
\begin{equation*}
T_{\Omega }(t)f=\left[ \cos \left( \Omega \frac{\xi _{3}}{|\xi |}t\right)
e^{-|\xi |^{2}t}I\widehat{f}(\xi )+\sin \left( \Omega \frac{\xi _{3}}{|\xi |}%
t\right) e^{-|\xi |^{2}t}\mathcal{R}(\xi )\widehat{f}(\xi )\right]^{\vee}
\end{equation*}%
for divergence-free vector fields $f$, where $I$ is the identity matrix in $%
\mathbb{R}^{3}$ and $\mathcal{R}(\xi )$ is the skew-symmetric matrix symbol
\begin{equation*}
\mathcal{R}(\xi )=\frac{1}{|\xi |}\left(
\begin{array}{ccc}
0 & \xi _{3} & -\xi _{2} \\
-\xi _{3} & 0 & \xi _{1} \\
\xi _{2} & -\xi _{1} & 0%
\end{array}%
\right) \ \text{ for }\ \xi \in \mathbb{R}^{3}\setminus \{0\}.
\end{equation*}%
Vector-fields $u$ satisfying the formulation (\ref{integralform}) are called
mild solutions for (\ref{NSC}).

In the last decades, the global well-posedness of models in fluid mechanics
has been studied by several authors of the mathematical community,
particularly in physical models of rotating fluids as the system (\ref{NSC}%
). In what follows, we give a brief review on some of these results. We
start with the works of Babin, Mahalov and Nicolaenko \cite%
{BabinMahalovNicolaenko1997,BabinMahalovNicolaenko1999,BabinMahalovNicolaenko2001}%
, who showed the global existence and regularity of solutions for (\ref{NSC}%
) with periodic initial velocity provided that the speed of rotation $%
|\Omega |$ is sufficiently large. In \cite%
{Chemin-Gallagher2002,Chemin-Gallagher2006}, Chemin\textit{\ et al}.
obtained a unique global strong Leray-type solution for large $|\Omega |$
and initial data $u_{0}(x)\in L^{2}(\mathbb{R}^{2})^{3}+H^{\frac{1}{2}}(%
\mathbb{R}^{3})^{3}$ (notice that the first parcel of $u_{0}(x)$ depends on $%
(x_{1},x_{2})$ where $x=(x_{1},x_{2},x_{3})$). For almost periodic initial
data and using the $l^{1}$-norm of amplitudes with sum closed frequency set,
Yoneda \cite{Yoneda2011} proved the existence of solutions for large times
and sufficiently large $|\Omega |$. Considering the mild (semigroup)
formulation, the global well-posedness in homogeneous Sobolev spaces $\dot{H}%
^{s}(\mathbb{R}^{3})$ with $1/2\leq s<3/4$ was obtained by Iwabuchi and
Takada \cite{Takada2013}. They considered sufficiently large $|\Omega |$
(depending on the size of $\Vert u_{0}\Vert _{\dot{H}^{s}}$) when $1/2<s<3/4$%
. In the critical case $s=1/2$, they used a class of precompact subsets in $%
\dot{H}^{1/2}(\mathbb{R}^{3})$ in order to get similar results. Local
versions ($T$ large but finite) of the results in \cite{Takada2013} can be
found in \cite{IwabuchiTakada2015} for $1/2<s<5/4.$

Another type of results for (\ref{NSC}) is the uniform global solvability
(or well-posedness) in which the smallness condition on $u_{0}$ is
independent of $|\Omega |$. Giga\textit{\ et al}. \cite{Giga2008} obtained
the uniform global solvability for small data $u_{0}$ in $FM_{0}^{-1}(%
\mathbb{R}^{3})=\mbox{div}(FM_{0}(\mathbb{R}^{3}))^{3}$, where $FM_{0}(%
\mathbb{R}^{3})$ denotes the space of the finite Radon measures with no
point mass at the origin. The space $FM_{0}^{-1}(\mathbb{R}^{3})$ is an
example of critical space for the 3D Navier-Stokes equations (NS) ((NSC)
with $\Omega =0$), i.e., its norm is invariant by the scaling $%
u_{0}^{\lambda }(x)\rightarrow \lambda u_{0}(\lambda x)$, for all $\lambda
>0 $. The uniform global well-posedness for small $u_{0}$ in the Sobolev
space $H^{\frac{1}{2}}(\mathbb{R}^{3})$ was proved by Hieber and Shibata
\cite{HieberShibata2010} and for small initial data in the critical
Fourier-Besov space $\dot{FB}_{p,\infty }^{2-\frac{3}{p}}(\mathbb{R}^{3})$
with $1<p\leq \infty $ and in $\dot{FB}_{1,1}^{-1}(\mathbb{R}^{3})\cap \dot{%
FB}_{1,1}^{0}(\mathbb{R}^{3})$ was proved by Konieczny and Yoneda \cite%
{Konieczny2011}. Iwabuchi and Takada \cite{IwabuchiTakada2014} obtained the
uniform global well-posedness with small initial velocity in the
Fourier-Besov $\dot{FB}_{1,2}^{-1}(\mathbb{R}^{3})$ as well as the
ill-posedness in $\dot{FB}_{1,q}^{-1}(\mathbb{R}^{3})$ for $2<q\leq \infty$.
These results were extended to the framework of critical
Fourier-Besov-Morrey spaces by Almeida, Ferreira and Lima \cite%
{Almeida-Fer-Lima}.

Concerning the asymptotic behavior for (\ref{NSC}), we quote the work of
Iwabuchi, Mahalov and Takada \cite{IwabuchiMahalovTakada2016}, where they
treated the high-rotating cases and proved the asymptotic stability of large
time periodic solutions for large initial perturbations. We also mention
\cite{Chemin-Gallagher2006} where the reader can find convergence results of
solutions towards a two-dimensional model as $|\Omega |\rightarrow \infty $
(see also references therein).

It is worthy to highlight that global existence of strong, mild or smooth
solutions for the Navier-Stokes equations ($\Omega =0$), without assume
smallness conditions on $u_{0}$, are outstanding open problems. Thus, global
solvability results for (\ref{NSC}) with arbitrary data in suitable spaces
show an interesting \textquotedblleft smoothing effect\textquotedblright\
due to the Coriolis parameter $\Omega $.

In this paper, we show the global well-posedness of (\ref{NSC}) for large $%
|\Omega |$ and arbitrary initial data $u_{0}$ belonging to homogeneous Besov
spaces $\dot{B}_{2,q}^{s}(\mathbb{R}^{3})$ where $1\leq q\leq \infty $ and $%
1/2\leq s<3/4$. In fact, for the cases $s\in (1/2,3/4)$ with $q=\infty $ and
$s=1/2$ with $q\in \lbrack 2,\infty ],$ we introduce the suitable
initial-data classes $\mathcal{I}$ and $\mathcal{F}_{0}$ (see (\ref%
{aux-space-77}) and (\ref{aux-space-772})), respectively, whose definitions
depend on the Stokes-Coriolis semigroup and Besov spaces. Also, we analyze
the asymptotic behavior of solutions as $\left\vert \Omega \right\vert
\rightarrow \infty $. For the case $1/2<s<3/4$, we use some space-time
estimates of Strichartz type for the Stokes-Coriolis semigroup, and also the
condition of $|\Omega |$ being large with respect to the $\dot{B}_{2,q}^{s}$%
-norm ($\mathcal{I}$-norm for $q=\infty $) of the initial data $u_{0}$ (a
power-type dependence). For the critical case $s=1/2$, $|\Omega |$ depends
on initial data belonging to precompact sets $D\subset \mathcal{F}_{0}$. In
view of the strict continuous inclusions $\dot{H}^{1/2}\subset \mathcal{F}%
_{0}$ and
\begin{equation*}
\dot{B}_{2,1}^{s}\subset \dot{B}_{2,q_{1}}^{s}\subset \dot{H}^{s}=\dot{B}%
_{2,2}^{s}\subset \dot{B}_{2,q_{2}}^{s}\subset \dot{B}_{2,\infty }^{s},
\end{equation*}%
for $1\leq q_{1}\leq 2\leq q_{2}<\infty ,$ our results provide a new initial
data class for the global well-posedness of (\ref{NSC}) and, in particular,
a class larger than that of \cite{Takada2013}.

Throughout this paper, we denote by $C>0$ constants that may differ even on
the same line. Also, the notation $C=C(a_{1},\ldots ,a_{k})$ indicates that $%
C$ depends on the quantities $a_{1},...,a_{k}.$

The outline of this paper is as follows. Section 2 is devoted to review some
basic facts about homogeneous Besov spaces and certain mixed space-time
functional settings. Estimates in Besov norms for the semigroup $T_{\Omega }$
and the Duhamel integral term in (\ref{integralform}) are the subject of
Section 3. In Section 4, we state and prove our global well-posedness and
asymptotic behavior results for (\ref{NSC}).

\section{Function spaces}

This section is devoted to some preliminaries about homogeneous Besov spaces
and some mixed space-time functional settings.

We start with the definition of the homogeneous Besov spaces. For this, let $%
\mathcal{S}(\mathbb{R}^{3})$ and $\mathcal{S}^{\prime }(\mathbb{R}^{3})$
stand for the Schwartz class and the space of tempered distributions,
respectively. Let $\widehat{f}$ denote the Fourier transform of $f\in
\mathcal{S}^{\prime }(\mathbb{R}^{3})$.

Consider a nonnegative radial function $\phi _{0}\in \mathcal{S}(\mathbb{R}%
^{3})$ satisfying
\begin{equation*}
0\leq \widehat{\phi }_{0}(\xi )\leq 1\text{ for all }\xi \in \mathbb{R}%
^{3},\ \mbox{supp}\ \widehat{\phi }_{0}\subset \{\xi \in \mathbb{R}^{3}:%
\frac{1}{2}\leq |\xi |\leq 2\}\text{ and }\sum_{j\in \mathbb{Z}}\widehat{%
\phi }_{j}(\xi )=1\ \ \mbox{for all}\ \ \xi \in \mathbb{R}^{3}\backslash
\{0\},
\end{equation*}%
where $\phi _{j}(x)=2^{3j}\phi _{0}(2^{j}x).$ For $f\in \mathcal{S}^{\prime
}(\mathbb{R}^{3}),$ the Littlewood-Paley operator $\{\Delta _{j}\}_{j\in
\mathbb{Z}}$ is defined by $\Delta _{j}f=\phi _{j}\ast f$.

Let $s\in \mathbb{R}$ and $1\leq p,q\leq \infty $ and let $\mathcal{P}$
denote the set of polynomials with $3$ variables. The homogeneous Besov
space, denoted by ${\dot{B}}_{p,q}^{s}(\mathbb{R}^{3})$, is defined as the
set of all $f\in \mathcal{S}^{\prime }(\mathbb{R}^{3})/\mathcal{P}$ such
that the following norm is finite
\begin{equation*}
\Vert f\Vert _{{\dot{B}}_{p,q}^{s}}=\Vert \{2^{sj}\Vert \Delta _{j}f\Vert
_{L^{p}}\}_{j\in \mathbb{Z}}\Vert _{l^{q}(\mathbb{Z})}.
\end{equation*}%
The pair $(\dot{B}_{p,q}^{s},\Vert \cdot \Vert _{{\dot{B}}_{p,q}^{s}})$ is a
Banach space. We will denote abusively distributions in $\mathcal{S}^{\prime
}(\mathbb{R}^{3})$ and their equivalence classes in $\mathcal{S}^{\prime }(%
\mathbb{R}^{3})/\mathcal{P}$ in the same way. The space $\mathcal{S}_{0}(%
\mathbb{R}^{3})$ of functions in $\mathcal{S}(\mathbb{R}^{3})$ whose Fourier
transforms are supported away from $0$ is dense in $\dot{B}_{p,q}^{s}(%
\mathbb{R}^{3})$ for $1\leq p,q<\infty $. For more details, see \cite%
{BahouriCheminDanchin2011}.

Using a duality argument, the norm $\Vert u\Vert _{\dot{B}_{p,q}^{s}}$ can
be estimated as follows
\begin{equation}
\Vert u\Vert _{\dot{B}_{p,q}^{s}}\leq C\sup_{\phi \in Q_{p^{\prime
},q^{\prime }}^{-s}}\left\vert \langle u,\phi \rangle \right\vert
\label{duality-1}
\end{equation}%
where $Q_{p^{\prime },q^{\prime }}^{-s}(\mathbb{R}^{3})$ denotes the set of
all functions $\phi \in \mathcal{S}(\mathbb{R}^{3})\cap \dot{B}_{p^{\prime
},q^{\prime }}^{-s}$ such that $\left\Vert \phi \right\Vert _{\dot{B}%
_{p^{\prime },q^{\prime }}^{-s}}\leq 1$ and $\langle \cdot ,\cdot \rangle $
is defined by
\begin{equation*}
\langle u,\phi \rangle :=\sum_{|j-j^{\prime }|\leq 1}\int_{\mathbb{R}%
^{3}}\Delta _{j}u(x)\Delta _{j^{\prime }}\phi (x)\ dx
\end{equation*}%
for $u\in \dot{B}_{p,q}^{s}(\mathbb{R}^{3})$ and $\phi \in Q_{p^{\prime
},q^{\prime }}^{-s}(\mathbb{R}^{3})$.

The next lemma contains a Leibniz type rule in the framework of Besov spaces.

\begin{lemma}[see \protect\cite{Chae2004}]
\label{product} Let $s>0,$ $1\leq q\leq \infty $ and $1\leq
p,p_{1},p_{2},r_{1},r_{2}\leq \infty $ be such that $\frac{1}{p}=\frac{1}{%
p_{1}}+\frac{1}{p_{2}}=\frac{1}{r_{1}}+\frac{1}{r_{2}}$. Then, there exists
a universal constant $C>0$ such that
\begin{equation*}
\Vert fg\Vert _{\dot{B}_{p,q}^{s}}\leq C\left( \Vert f\Vert
_{L^{p_{1}}}\Vert g\Vert _{\dot{B}_{p_{2},q}^{s}}+\Vert g\Vert
_{L^{r_{1}}}\Vert f\Vert _{\dot{B}_{r_{2},q}^{s}}\right) .
\end{equation*}
\end{lemma}

Considering in particular $p=r$ and $p_{i}=r_{i}$ in Lemma \ref{product}, we
have that
\begin{equation*}
\Vert fg\Vert _{\dot{B}_{r,q}^{s}}\leq C\left( \Vert f\Vert
_{L^{r_{1}}}\Vert g\Vert _{\dot{B}_{r_{2},q}^{s}}+\Vert g\Vert
_{L^{r_{1}}}\Vert f\Vert _{\dot{B}_{r_{2},q}^{s}}\right) .
\end{equation*}%
If $\frac{1}{r}=\frac{2}{r_{2}}-\frac{s}{3}$, then $\frac{1}{r_{1}}=\frac{2}{%
r_{2}}-\frac{s}{3}$ and we can use the embedding $\dot{B}_{r_{2},q}^{s}(%
\mathbb{R}^{3})\hookrightarrow L^{r_{1}}(\mathbb{R}^{3})$ to obtain
\begin{equation}
\Vert fg\Vert _{\dot{B}_{r,q}^{s}}\leq C\Vert f\Vert _{\dot{B}%
_{r_{2},q}^{s}}\Vert g\Vert _{\dot{B}_{r_{2},q}^{s}}.  \label{remark1}
\end{equation}%
The reader is referred to \cite{BergLofstrom1976} for more details on $\dot{B%
}_{p,q}^{s}$-spaces and their properties.\bigskip

We finish this section by recalling some mixed space-time functional spaces.
Let $\theta \geq 1$, we denote by $L^{\theta }(0,\infty ;\dot{B}_{p,q}^{s}(%
\mathbb{R}^{3}))$ the set of all distributions $f$ such that
\begin{equation*}
\Vert f\Vert _{L^{\theta }(0,\infty ;\dot{B}_{p,q}^{s})}=\left\Vert \Vert
f(t)\Vert _{\dot{B}_{p,q}^{s}}\right\Vert _{L_{t}^{\theta }(0,\infty
)}<\infty .
\end{equation*}%
Also, we denote by $\tilde{L}^{\theta }(0,\infty ;\dot{B}_{p,q}^{s}(\mathbb{R%
}^{3}))$ the set of all distributions $f$ such that
\begin{equation*}
\Vert f\Vert _{\tilde{L}^{\theta }(0,\infty ;\dot{B}_{p,q}^{s})}=\left\Vert
\{2^{js}\Vert \Delta _{j}f\Vert _{L^{\theta }(0,\infty ;L^{p})}\}_{j\in
\mathbb{Z}}\right\Vert _{l^{q}(\mathbb{Z})}<\infty .
\end{equation*}%
As consequence of the Minkowski inequality, we have the following embeddings
\begin{equation}
\begin{split}
& L^{\theta }(0,\infty ;\dot{B}_{p,q}^{s})\hookrightarrow \tilde{L}^{\theta
}(0,\infty ;\dot{B}_{p,q}^{s}),\text{ if }\theta \leq q, \\
& \tilde{L}^{\theta }(0,\infty ;\dot{B}_{p,q}^{s})\hookrightarrow L^{\theta
}(0,\infty ;\dot{B}_{p,q}^{s}),\text{ if }\theta \geq q.
\end{split}
\label{embedding}
\end{equation}

\section{\protect\bigskip Estimates}

Firstly, we recall some estimates for the heat semigroup $e^{t\Delta }$ in
Besov spaces \cite{KozonoOgawaTanuichi2003} and the dispersive estimates for
$T_{\Omega }(t)$ obtained in \cite{IwabuchiTakada2013D}.

\begin{lemma}[see \protect\cite{KozonoOgawaTanuichi2003}]
\label{heatbesov} Let $-\infty <s_{0}\leq s_{1}<\infty $, $1\leq p,q\leq
\infty $ and $f\in \dot{B}_{p,q}^{s_{0}}(\mathbb{R}^{3})$. Then, there
exists a positive constant $C=C(s_{0},s_{1})$ such that
\begin{equation*}
\Vert e^{t\Delta }f\Vert _{\dot{B}_{p,q}^{s_{1}}}\leq Ct^{-\frac{3}{2}%
(s_{1}-s_{0})}\Vert f\Vert _{\dot{B}_{p,q}^{s_{0}}}\text{, for all }t>0.
\end{equation*}
\end{lemma}

Before stating the dispersive estimates of \cite{IwabuchiTakada2013D}, we
need to define the operators
\begin{equation}
\mathcal{G}_{\pm }(\tau )[f]=\left[ e^{\pm i\tau \frac{\xi _{3}}{|\xi |}}%
\widehat{f}\right]^{\vee} ,\text{ for }\tau \in \mathbb{R}\text{,}
\label{definitionoperatorG}
\end{equation}%
and the matrix $\mathcal{R}$ of singular integral operators
\begin{equation}
\mathcal{R}=\left(
\begin{array}{ccc}
0 & R_{3} & -R_{2} \\
-R_{3} & 0 & R_{1} \\
R_{2} & -R_{1} & 0%
\end{array}%
\right) .  \label{matrix-1}
\end{equation}%
Using (\ref{definitionoperatorG}) and (\ref{matrix-1}), $T_{\Omega }(t)$ can
be expressed as%
\begin{equation}
T_{\Omega }(t)f=\frac{1}{2}\mathcal{G}_{+}(\Omega t)[e^{t\Delta }(I+\mathcal{%
R})f]+\frac{1}{2}\mathcal{G}_{-}(\Omega t)[e^{t\Delta }(I-\mathcal{R})f]
\label{semigroup-formula-1}
\end{equation}%
for $t\geq 0$ and $\Omega \in \mathbb{R}$. \

Notice that the operators $\mathcal{G}_{\pm }(t\Omega )$ correspond to the
oscillating parts of $T_{\Omega }(t)$.

\begin{lemma}[see \protect\cite{IwabuchiTakada2013D}]
\label{operatorG} Let $s,t\in \mathbb{R}$, $2\leq p\leq \infty $ and $f\in
\dot{B}_{p^{\prime },q}^{s+3\left( 1-\frac{2}{p}\right) }(\mathbb{R}^{3})$
with $\frac{1}{p}+\frac{1}{p^{\prime }}=1$. Then, there exists a constant $%
C=C(p)>0$ such that%
\begin{equation*}
\Vert \mathcal{G}_{\pm }(t)[f]\Vert _{\dot{B}_{p,q}^{s}}\leq C\left( \frac{%
\log {(e+|t|)}}{1+|t|}\right) ^{\frac{1}{2}\left( 1-\frac{2}{p}\right)
}\Vert f\Vert _{\dot{B}_{p^{\prime },q}^{s+3\left( 1-\frac{2}{p}\right) }}.
\end{equation*}%
\bigskip
\end{lemma}

In what follows, we establish our estimates in Besov spaces for $T_{\Omega
}(t)$ and the Duhamel term $\int_{0}^{t}T_{\Omega }(t-\tau )\mathbb{P}\nabla
f(\tau )\ d\tau $. We start with three lemmas for $T_{\Omega }(t).$

\begin{lemma}
\label{linearestimative1} Assume that $s,\Omega \in \mathbb{R}$, $t>0,$ $%
1<r\leq p^{\prime }\leq 2\leq p<\infty $ and $1\leq q\leq \infty $, and let $%
k$ be a multi-index. Then, there exists a constant $C>0$ (independent of $%
\Omega $ and $t$ ) such that
\begin{equation*}
\Vert \nabla _{x}^{k}T_{\Omega }(t)f\Vert _{\dot{B}_{p,q}^{s}}\leq C\left(
\frac{\log {\ (e+|t\Omega |)}}{1+|t\Omega |}\right) ^{\frac{1}{2}\left( 1-%
\frac{2}{p}\right) }t^{-\frac{\left\vert k\right\vert }{2}-\frac{3}{2}\left(
\frac{1}{r}-\frac{1}{p}\right) }\Vert f\Vert _{\dot{B}_{r,q}^{s}},
\end{equation*}%
for all $f\in \dot{B}_{r,q}^{s}(\mathbb{R}^{3})$.
\end{lemma}

\textbf{Proof. }Using the representation (\ref{semigroup-formula-1}), Lemma %
\ref{operatorG}, the embedding $\dot{B}_{r,q}^{s+3\left( \frac{1}{r}-\frac{1%
}{p}\right) }(\mathbb{R}^{3})\hookrightarrow \dot{B}_{p^{\prime
},q}^{s+3\left( 1-\frac{2}{p}\right) }(\mathbb{R}^{3})$ and Lemma \ref%
{heatbesov}, we obtain

\begin{equation*}
\begin{split}
\Vert \nabla _{x}^{k}T_{\Omega }(t)f\Vert _{\dot{B}_{p,q}^{s}}& \leq C\Vert
\mathcal{G}_{\pm }(t\Omega )[\nabla _{x}^{k}e^{t\Delta }f]\Vert _{\dot{B}%
_{p,q}^{s}} \\
& \leq C\left( \frac{\log {(e+|t\Omega |)}}{1+|t\Omega |}\right) ^{\frac{1}{2%
}\left( 1-\frac{2}{p}\right) }\Vert \nabla _{x}^{k}e^{t\Delta }f\Vert _{\dot{%
B}_{p^{\prime },q}^{s+3\left( 1-\frac{2}{p}\right) }} \\
& \leq C\left( \frac{\log {(e+|t\Omega |)}}{1+|t\Omega |}\right) ^{\frac{1}{2%
}\left( 1-\frac{2}{p}\right) }\Vert \nabla _{x}^{k}e^{t\Delta }f\Vert _{\dot{%
B}_{r,q}^{s+3\left( \frac{1}{r}-\frac{1}{p}\right) }} \\
& \leq C\left( \frac{\log {(e+|t\Omega |)}}{1+|t\Omega |}\right) ^{\frac{1}{2%
}\left( 1-\frac{2}{p}\right) }t^{-\frac{\left\vert k\right\vert }{2}-\frac{3%
}{2}\left( \frac{1}{r}-\frac{1}{p}\right) }\Vert f\Vert _{\dot{B}_{r,q}^{s}}.
\end{split}%
\end{equation*}%
\fin

\begin{lemma}
Let $1\leq q<\infty $. Consider $s,p,\theta \in \mathbb{R}$ satisfying
\begin{equation*}
0\leq s<\frac{3}{p},\ \ 2<p<6\ \ \text{ and }\ \ \frac{3}{4}-\frac{3}{2p}%
\leq \frac{1}{\theta }<\min \left\{ \frac{1}{2},1-\frac{2}{p},\frac{1}{q}%
\right\} .
\end{equation*}%
Then, there exists $C>0$ (independent of $t\geq 0$ and $\Omega \in \mathbb{R}
$) such that

\begin{equation}
\Vert T_{\Omega }(t)f\Vert _{L^{\theta }(0,\infty ;\dot{B}_{p,q}^{s})}\leq
C|\Omega |^{-\frac{1}{\theta }+\frac{3}{4}-\frac{3}{2p}}\Vert f\Vert _{\dot{B%
}_{2,q}^{s}},  \label{linearestimative2}
\end{equation}%
for all $f\in \dot{B}_{2,q}^{s}(\mathbb{R}^{3})$.
\end{lemma}

\textbf{Proof. } By duality and estimate (\ref{duality-1}), notice that (\ref%
{linearestimative2}) holds true provided that
\begin{equation}
\begin{split}
I& :=\left\vert \int_{0}^{\infty }\sum_{|j-k|\leq 1}\int_{\mathbb{R}%
^{3}}\Delta _{j}\mathcal{G}_{\pm }(\Omega t)[e^{t\Delta }f](x)\overline{%
\Delta _{k}\phi (x,t)}\ dxdt\right\vert \\
& \leq C|\Omega |^{-\frac{1}{\theta }+\frac{3}{4}-\frac{3}{2p}}\Vert f\Vert
_{\dot{B}_{2,q}^{s}}\Vert \phi \Vert _{L^{\theta ^{\prime }}(0,\infty ;\dot{B%
}_{p^{\prime },q^{\prime }}^{-s})},
\end{split}
\label{est-aux-2}
\end{equation}%
for all $\phi \in C_{0}^{\infty }(\mathbb{R}^{3}\times (0,\infty ))$ with $%
0\notin $supp$(\widehat{\phi }(\xi ,t))$ for each $t>0,$ where $%
1/p+1/p^{\prime }=1$, $1/\theta +1/\theta ^{\prime }=1$ and $1/q+1/q^{\prime
}=1.$

For (\ref{est-aux-2}), we use Parseval formula, H\"{o}lder inequality, the
inclusion $\dot{B}_{p,2}^{0}(\mathbb{R}^{3})\hookrightarrow L^{p}(\mathbb{R}%
^{3})$ and Lemma \ref{linearestimative1} in order to estimate%
\begin{equation}
\begin{split}
I& \leq \sum_{|j-k|\leq 1}\left\vert \int_{0}^{\infty }\int_{\mathbb{R}%
^{3}}\Delta _{j}\mathcal{G}_{\pm }(\Omega t)[e^{t\Delta }f](x)\overline{%
\Delta _{k}\phi (x,t)}\ dxdt\right\vert \\
& =\sum_{|j-k|\leq 1}\left\vert \int_{0}^{\infty }\int_{\mathbb{R}%
^{3}}\Delta _{j}f(x)\overline{\Delta _{k}\mathcal{G}_{\mp }(\Omega
t)[e^{t\Delta }\phi (t)](x)}\ dxdt\right\vert \\
& =\sum_{|j-k|\leq 1}\left\vert \int_{\mathbb{R}^{3}}\Delta
_{j}f(x)\int_{0}^{\infty }\overline{\Delta _{k}\mathcal{G}_{\mp }(\Omega
t)[e^{t\Delta }\phi (t)](x)}\ dtdx\right\vert \\
& \leq C\sum_{|j-k|\leq 1}\Vert \Delta _{j}f\Vert _{L^{2}}\left\Vert
\int_{0}^{\infty }\Delta _{k}\mathcal{G}_{\mp }(\Omega t)[e^{t\Delta }\phi
(t)]\ dt\right\Vert _{L^{2}} \\
& \leq C2^{|s|}\sum_{|j-k|\leq 1}2^{js}\Vert \Delta _{j}f\Vert
_{L^{2}}2^{-ks}\left\Vert \int_{0}^{\infty }\Delta _{k}\mathcal{G}_{\mp
}(\Omega t)[e^{t\Delta }\phi (t)]\ dt\right\Vert _{L^{2}} \\
& \leq C2^{|s|}\Vert f\Vert _{\dot{B}_{2,q}^{s}}\left( \sum_{k\in \mathbb{Z}%
}2^{-ksq^{\prime }}\left\Vert \int_{0}^{\infty }\Delta _{k}\mathcal{G}_{\mp
}(\Omega t)[e^{t\Delta }\phi (t)]\ dt\right\Vert _{L^{2}}^{q^{\prime
}}\right) ^{\frac{1}{q^{\prime }}}.
\end{split}
\label{12}
\end{equation}

Now, we are going to prove that
\begin{equation*}
I_{k}^{2}\leq C|\Omega |^{-\frac{1}{\theta }+\frac{3}{4}-\frac{3}{2p}}\Vert
\Delta _{k}\phi \Vert _{L^{\theta ^{\prime }}(0,\infty ;L^{p^{\prime
}})}^{2},
\end{equation*}%
where
\begin{equation*}
I_{k}:=\left\Vert \int_{0}^{\infty }\Delta _{k}\mathcal{G}_{\mp }(\Omega
t)[e^{t\Delta }\phi (t)]\ dt\right\Vert _{L^{2}}.
\end{equation*}%
In fact, using the Parseval formula, H\"{o}lder inequality, the embedding $%
\dot{B}_{p,2}^{0}(\mathbb{R}^{3})\hookrightarrow L^{p}(\mathbb{R}^{3})$ and
Lemma \ref{linearestimative1}, we have
\begin{equation*}
\begin{split}
I_{k}^{2}& =\left\langle \int_{0}^{\infty }\Delta _{k}\mathcal{G}_{\mp
}(\Omega t)[e^{t\Delta }\phi (t)]\ dt,\int_{0}^{\infty }\Delta _{k}\mathcal{G%
}_{\mp }(\Omega \tau )[e^{\tau \Delta }\phi (\tau )]\ d\tau \right\rangle
_{L^{2}} \\
& =\int_{0}^{\infty }\int_{0}^{\infty }\int_{\mathbb{R}^{3}}\Delta _{k}%
\mathcal{G}_{\mp }(\Omega t)[e^{t\Delta }\phi (t)](x)\overline{\Delta _{k}%
\mathcal{G}_{\mp }(\Omega \tau )[e^{\tau \Delta }\phi (\tau )](x)}\ dxd\tau
dt \\
& \leq \int_{0}^{\infty }\int_{0}^{\infty }\Vert \Delta _{k}\phi (t)\Vert
_{L^{p^{\prime }}}\Vert \Delta _{k}\mathcal{G}_{\pm }(\Omega (t-\tau
))[e^{(t+\tau )\Delta }\phi (\tau )]\Vert _{L^{p}}\ d\tau dt \\
& \leq C\int_{0}^{\infty }\int_{0}^{\infty }\Vert \Delta _{k}\phi (t)\Vert
_{L^{p^{\prime }}}\Vert \Delta _{k}\mathcal{G}_{\pm }(\Omega (t-\tau
))[e^{(t+\tau )\Delta }\phi (\tau )]\Vert _{\dot{B}_{p,2}^{0}}\ d\tau dt \\
& \leq C\int_{0}^{\infty }\int_{0}^{\infty }\Vert \Delta _{k}\phi (t)\Vert
_{L^{p^{\prime }}}\left( \frac{\log {(e+|\Omega ||t-\tau |)}}{1+|\Omega
||t-\tau |}\right) ^{\frac{1}{2}\left( 1-\frac{2}{p}\right) }\Vert
e^{(t+\tau )\Delta }\Delta _{k}\phi (\tau )\Vert _{\dot{B}_{p^{\prime
},2}^{3\left( 1-\frac{2}{p}\right) }}\ d\tau dt.
\end{split}%
\end{equation*}%
By Lemma \ref{heatbesov} and the embedding $L^{p^{\prime }}(\mathbb{R}%
^{3})\hookrightarrow \dot{B}_{p^{\prime },2}^{0}(\mathbb{R}^{3})$ for $%
p^{\prime }<2$, it follows that
\begin{equation*}
\begin{split}
\Vert e^{(t+\tau )\Delta }\Delta _{k}\phi (\tau )\Vert _{\dot{B}_{p^{\prime
},2}^{3\left( 1-\frac{2}{p}\right) }}& \leq C(t+\tau )^{-\frac{3}{2}\left( 1-%
\frac{2}{p}\right) }\Vert \Delta _{k}\phi (\tau )\Vert _{\dot{B}_{p^{\prime
},2}^{0}} \\
& \leq C|t-\tau |^{-\frac{3}{2}\left( 1-\frac{2}{p}\right) }\Vert \Delta
_{k}\phi (\tau )\Vert _{L^{p^{\prime }}}.
\end{split}%
\end{equation*}%
Thus,
\begin{equation}
\begin{split}
I_{k}^{2}& \leq C\int_{0}^{\infty }\int_{0}^{\infty }\Vert \Delta _{k}\phi
(t)\Vert _{L^{p^{\prime }}}\left( \frac{\log {(e+|\Omega ||t-\tau |)}}{%
1+|\Omega ||t-\tau |}\right) ^{\frac{1}{2}\left( 1-\frac{2}{p}\right)
}|t-\tau |^{-\frac{3}{2}\left( 1-\frac{2}{p}\right) }\Vert \Delta _{k}\phi
(\tau )\Vert _{L^{p^{\prime }}}\ d\tau dt \\
& \leq C\Vert \Delta _{k}\phi \Vert _{L^{\theta ^{\prime }}(0,\infty
;L^{p^{\prime }})}\left\Vert \int_{0}^{\infty }h(\cdot -\tau )\Vert \Delta
_{k}\phi (\tau )\Vert _{L^{p^{\prime }}}d\tau \right\Vert _{L^{\theta
}(0,\infty )},
\end{split}
\label{11}
\end{equation}%
where
\begin{equation*}
h(t)=\left( \frac{\log {(e+|\Omega ||t|)}}{1+|\Omega ||t|}\right) ^{\frac{1}{%
2}\left( 1-\frac{2}{p}\right) }|t|^{-\frac{3}{2}\left( 1-\frac{2}{p}\right)
}.
\end{equation*}%
We consider the cases $\frac{1}{\theta }>\frac{3}{4}-\frac{3}{2p}$ and $%
\frac{1}{\theta }=\frac{3}{4}-\frac{3}{2p}$. In the first case, notice that
\begin{equation*}
\Vert h\Vert _{L^{\frac{\theta }{2}}}=C|\Omega |^{-\frac{2}{\theta }+\frac{3%
}{2}-\frac{3}{p}}.
\end{equation*}%
Therefore, using Young inequality in (\ref{11}) and the above equality, we
obtain
\begin{equation*}
I_{k}^{2}\leq C|\Omega |^{-\frac{2}{\theta }+\frac{3}{2}-\frac{3}{p}}\Vert
\Delta _{k}\phi \Vert _{L^{\theta ^{\prime }}(0,\infty ;L^{p^{\prime
}})}^{2}.
\end{equation*}%
Now, multiplying by $2^{-ks}$, applying the $l^{q^{\prime }}(\mathbb{Z})$%
-norm and using (\ref{embedding}), we arrive at
\begin{equation}
\begin{split}
\left( \sum_{k\in \mathbb{Z}}2^{-ksq^{\prime }}I_{k}^{q^{\prime }}\right) ^{%
\frac{1}{q^{\prime }}}& \leq C|\Omega |^{-\frac{1}{\theta }+\frac{3}{4}-%
\frac{3}{2p}}\left( \sum_{k\in \mathbb{Z}}2^{-ksq^{\prime }}\Vert \Delta
_{k}\phi \Vert _{L^{\theta ^{\prime }}(0,\infty ;L^{p^{\prime
}})}^{q^{\prime }}\right) ^{\frac{1}{q^{\prime }}} \\
& \leq C|\Omega |^{-\frac{1}{\theta }+\frac{3}{4}-\frac{3}{2p}}\Vert \phi
\Vert _{L^{\theta ^{\prime }}(0,\infty ;\dot{B}_{p^{\prime },q^{\prime
}}^{-s})}.
\end{split}
\label{13}
\end{equation}%
It follows from (\ref{12}) and (\ref{13}) that
\begin{equation}
I\leq C|\Omega |^{-\frac{1}{\theta }+\frac{3}{4}-\frac{3}{2p}}\Vert f\Vert _{%
\dot{B}_{r,q}^{s}}\Vert \phi \Vert _{L^{\theta ^{\prime }}(0,\infty ;\dot{B}%
_{p^{\prime },q^{\prime }}^{-s})},  \label{15a}
\end{equation}%
with $C>0$ independent of $\phi $ and $f$.

In the second case $\frac{1}{\theta }=\frac{3}{4}-\frac{3}{2p}$, we use the
fact $h(t)\leq |t|^{-\frac{3}{2}\left( 1-\frac{2}{p}\right) }$ and
Hardy-Littlewood-Sobolev inequality in (\ref{11}) to obtain
\begin{equation}
I_{k}^{2}\leq C\Vert \Delta _{k}\phi \Vert _{L^{\theta ^{\prime }}(0,\infty
;L^{p^{\prime }})}^{2}.  \label{14}
\end{equation}%
Thus, using (\ref{14}) and proceeding as in (\ref{13}), we obtain a constant
$C>0$ (independent of $\phi $ and $f$) such that%
\begin{equation}
I\leq C\Vert f\Vert _{\dot{B}_{2,q}^{s}}\Vert \phi \Vert _{L^{\theta
^{\prime }}(0,\infty ;\dot{B}_{p^{\prime },q^{\prime }}^{-s})}.  \label{15b}
\end{equation}%
Estimates (\ref{15a}) and (\ref{15b}) give the desired result. \fin

\begin{lemma}
\label{critical_linear_semigroup} Assume that $1\leq q<4 $ and $f\in \dot{B}%
_{2,q}^{\frac{1}{2}}(\mathbb{R}^{3})$. Then,%
\begin{equation}
\lim_{|\Omega |\rightarrow \infty }\Vert T_{\Omega }(\cdot )f\Vert
_{L^{4}(0,\infty ;\dot{B}_{3,q}^{\frac{1}{2}})}=0.  \label{3.12}
\end{equation}
\end{lemma}

\textbf{Proof. } Since $\overline{\mathcal{S}_{0}(\mathbb{R}^{3})}%
^{\left\Vert \cdot \right\Vert _{\dot{B}_{2,q}^{1/2}}}=\dot{B}_{2,q}^{\frac{1%
}{2}}$ for $q\neq \infty $ (see Section 2), there exists $(w_{k})_{k\in
\mathbb{N}}$ in $\mathcal{S}_{0}(\mathbb{R}^{3})$ such that $%
w_{k}\rightarrow f$ in $\dot{B}_{2,q}^{\frac{1}{2}}(\mathbb{R}^{3})$ as $%
k\rightarrow \infty $. Next, using Lemma \ref{linearestimative2}, we obtain

\begin{equation}
\begin{split}
\limsup_{\left\vert \Omega \right\vert \rightarrow \infty }\Vert T_{\Omega
}(\cdot )f\Vert _{L^{4}(0,\infty ;\dot{B}_{3,q}^{\frac{1}{2}})}& \leq
\limsup_{\left\vert \Omega \right\vert \rightarrow \infty }\Vert T_{\Omega
}(\cdot )(f-w_{k})\Vert _{L^{4}(0,\infty ;\dot{B}_{3,q}^{\frac{1}{2}%
})}+\limsup_{\left\vert \Omega \right\vert \rightarrow \infty }\Vert
T_{\Omega }(\cdot )w_{k}\Vert _{L^{4}(0,\infty ;\dot{B}_{3,q}^{\frac{1}{2}})}
\\
& \leq C\Vert w_{k}-f\Vert _{\dot{B}_{2,q}^{\frac{1}{2}}}+\limsup_{\left%
\vert \Omega \right\vert \rightarrow \infty }\Vert T_{\Omega }(\cdot
)w_{k}\Vert _{L^{4}(0,\infty ;\dot{B}_{3,q}^{\frac{1}{2}})}.
\end{split}
\label{3.13}
\end{equation}%
Choosing $p\in (\frac{8}{3},3),$ we have the conditions
\begin{equation*}
\frac{3}{4}-\frac{3}{2p}<\frac{1}{4}<\min \left\{ 1-\frac{2}{p},\frac{1}{q}%
\right\} \ \ \text{and}\ \ \frac{1}{2}-\frac{3}{2p}<0.
\end{equation*}%
Then, we can use $\dot{B}_{p,q}^{-\frac{1}{2}+\frac{3}{p}}(\mathbb{R}%
^{3})\hookrightarrow \dot{B}_{3,q}^{\frac{1}{2}}(\mathbb{R}^{3})$ and Lemma %
\ref{linearestimative2} to estimate

\begin{equation}
\begin{split}
\limsup_{\left\vert \Omega \right\vert \rightarrow \infty }\Vert T_{\Omega
}(\cdot )w_{k}\Vert _{L^{4}(0,\infty ;\dot{B}_{3,q}^{\frac{1}{2}})}& \leq
C\limsup_{\left\vert \Omega \right\vert \rightarrow \infty }\Vert T_{\Omega
}(\cdot )w_{k}\Vert _{L^{4}(0,\infty ;\dot{B}_{p,q}^{-\frac{1}{2}+\frac{3}{p}%
})} \\
& \leq C|\Omega |^{\frac{1}{2}-\frac{3}{2p}}\Vert w_{k}\Vert _{\dot{B}%
_{2,q}^{-\frac{1}{2}+\frac{3}{p}}}\rightarrow 0,\ \ \text{as}\ \ |\Omega
|\rightarrow \infty .
\end{split}
\label{3.15}
\end{equation}%
By (\ref{3.13}), (\ref{3.15}) and $\left\Vert w_{k}-f\right\Vert _{\dot{B}%
_{2,q}^{\frac{1}{2}}}\rightarrow 0$, it follows (\ref{3.12}). \fin\bigskip

The next two lemmas are concerned with the Duhamel term $\int_{0}^{t}T_{%
\Omega }(t-\tau )\mathbb{P}\nabla f(\tau )\ d\tau .$

\begin{lemma}
\label{nonlinear1} Let $s\in \mathbb{R}$ and $\Omega \in \mathbb{R}%
\backslash \{0\}$ and let $p,r,q,\theta $ be real numbers satisfying
\begin{gather*}
2<p<3,\ \ \ \frac{6}{5}<r<2,\ \ \ 1\leq q\leq \infty ,\ \ \ 1-\frac{1}{p}%
\leq \frac{1}{r}<\frac{1}{3}+\frac{1}{p}, \\
\max \left\{ 0,\frac{1}{2}-\frac{3}{2}\left( \frac{1}{r}-\frac{1}{p}\right) -%
\frac{1}{2}\left( 1-\frac{2}{p}\right) \right\} <\frac{1}{\theta }\leq \frac{%
1}{2}-\frac{3}{2}\left( \frac{1}{r}-\frac{1}{p}\right) .
\end{gather*}%
Then, there exists a universal constant $C>0$ such that
\begin{equation}
\left\Vert \int_{0}^{t}T_{\Omega }(t-\tau )\mathbb{P}\nabla f(\tau )\ d\tau
\right\Vert _{L^{\theta }(0,\infty ;\dot{B}_{p,q}^{s})}\leq C|\Omega |^{-%
\frac{1}{2}+\frac{3}{2}\left( \frac{1}{r}-\frac{1}{p}\right) +\frac{1}{%
\theta }}\Vert f\Vert _{L^{\frac{\theta }{2}}(0,\infty ;\dot{B}_{r,q}^{s})}.
\label{aux-lemma-nonlinear-1}
\end{equation}
\end{lemma}

\textbf{Proof. } Using Lemma \ref{linearestimative1} it follows that
\begin{equation}
\begin{split}
\Bigl\Vert \int_{0}^{t}T_{\Omega }(t-\tau )& \mathbb{P}\nabla f(\tau )\
d\tau \Bigr\Vert _{L^{\theta }(0,\infty ;\dot{B}_{p,q}^{s})}\leq C\left\Vert
\int_{0}^{t}\left\Vert T_{\Omega }(t-\tau )\mathbb{P}\nabla f(\tau
)\right\Vert _{\dot{B}_{p,q}^{s}}\ d\tau \right\Vert _{L^{\theta }(0,\infty
)} \\
& \leq C\left\Vert \int_{0}^{t}(t-\tau )^{-\frac{1}{2}-\frac{3}{2}\left(
\frac{1}{r}-\frac{1}{p}\right) }\left( \frac{\log {(e+|\Omega ||t-\tau |)}}{%
1+|\Omega ||t-\tau |}\right) ^{\frac{1}{2}\left( 1-\frac{2}{p}\right)
}\left\Vert f(\tau )\right\Vert _{\dot{B}_{r,q}^{s}}\ d\tau \right\Vert
_{L^{\theta }(0,\infty )}.
\end{split}
\label{16a}
\end{equation}%
We are going to prove (\ref{nonlinear1}) in two cases. First we consider the
case $\frac{1}{\theta }=\frac{1}{2}-\frac{3}{2}\left( \frac{1}{r}-\frac{1}{p}%
\right) $. Here, we note that
\begin{equation*}
\left( \frac{\log {(e+|\Omega ||t-\tau |)}}{1+|\Omega ||t-\tau |}\right) ^{%
\frac{1}{2}\left( 1-\frac{2}{p}\right) }\leq 1
\end{equation*}%
and employ Hardy-Littlewood-Sobolev inequality to estimate
\begin{equation}
\left\Vert \int_{0}^{t}(t-\tau )^{-\frac{1}{2}-\frac{3}{2}\left( \frac{1}{r}-%
\frac{1}{p}\right) }\left( \frac{\log {(e+|\Omega ||t-\tau |)}}{1+|\Omega
||t-\tau |}\right) ^{\frac{1}{2}\left( 1-\frac{2}{p}\right) }\left\Vert
f(\tau )\right\Vert _{\dot{B}_{r,q}^{s}}\ d\tau \right\Vert _{L^{\theta
}(0,\infty )}\leq C\Vert f\Vert _{L^{\frac{\theta }{2}}(0,\infty ;\dot{B}%
_{r,q}^{s})}.  \label{16}
\end{equation}%
Consider now the case $\frac{1}{\theta }<\frac{1}{2}-\frac{3}{2}\left( \frac{%
1}{r}-\frac{1}{p}\right) .$ Selecting $\ell $ such that $\frac{1}{\theta }=%
\frac{1}{\ell }+\frac{2}{\theta }-1$, a direct computation gives

\begin{equation}
\left\Vert (t-\tau )^{-\frac{1}{2}-\frac{3}{2}\left( \frac{1}{r}-\frac{1}{p}%
\right) }\left( \frac{\log {(e+|\Omega ||t-\tau |)}}{1+|\Omega ||t-\tau |}%
\right) ^{\frac{1}{2}\left( 1-\frac{2}{p}\right) }\right\Vert _{L^{\ell
}(0,\infty )}=C|\Omega |^{\frac{1}{\theta }-\frac{1}{2}+\frac{3}{2}\left(
\frac{1}{r}-\frac{1}{p}\right) }.  \label{aux-est-17}
\end{equation}%
By Young inequality and (\ref{aux-est-17}), we have that
\begin{equation}
\begin{split}
\Biggl\|& \int_{0}^{t}(t-\tau )^{-\frac{1}{2}-\frac{3}{2}\left( \frac{1}{r}-%
\frac{1}{p}\right) }\left( \frac{\log {(e+|\Omega ||t-\tau |)}}{1+|\Omega
||t-\tau |}\right) ^{\frac{1}{2}\left( 1-\frac{2}{p}\right) }\left\Vert
f(\tau )\right\Vert _{\dot{B}_{r,q}^{s}}\ d\tau \Biggr\|_{L^{\theta
}(0,\infty )} \\
& \leq \left\Vert t^{-\frac{1}{2}-\frac{3}{2}\left( \frac{1}{r}-\frac{1}{p}%
\right) }\left( \frac{\log {(e+|\Omega |t)}}{1+|\Omega |t}\right) ^{\frac{1}{%
2}\left( 1-\frac{2}{p}\right) }\right\Vert _{L^{\ell }(0,\infty )}\Vert
f\Vert _{L^{\frac{\theta }{2}}(0,\infty ;\dot{B}_{r,q}^{s})} \\
& =C|\Omega |^{\frac{1}{\theta }-\frac{1}{2}+\frac{3}{2}\left( \frac{1}{r}-%
\frac{1}{p}\right) }\Vert f\Vert _{L^{\frac{\theta }{2}}(0,\infty ;\dot{B}%
_{r,q}^{s})}.
\end{split}
\label{17}
\end{equation}%
The proof is completed by substituting (\ref{16}) and (\ref{17})\ into (\ref%
{16a}). \fin

\begin{lemma}
\label{nonlinearestimativecritical} Let $s,\Omega \in \mathbb{R}$ and $2\leq
q\leq \infty $. Then, there exists a universal constant $C>0$ such that
\begin{equation}
\left\Vert \int_{0}^{t}T_{\Omega }(t-\tau )\nabla f(\tau )\ d\tau
\right\Vert _{L^{\infty }(0,\infty ;\dot{B}_{2,q}^{s})\cap L^{4}(0,\infty ;%
\dot{B}_{3,q}^{s})}\leq C\Vert f\Vert _{L^{2}(0,\infty ;\dot{B}_{2,q}^{s})}.
\label{aux-lemma-10}
\end{equation}
\end{lemma}

\textbf{Proof. } We denote $X=X_{1}\cap X_{2}$ where $X_{1}=L^{\infty
}(0,\infty ;\dot{B}_{2,q}^{s})$ and $X_{2}=L^{4}(0,\infty ;\dot{B}%
_{3,q}^{s}) $. We start with estimates for the $X_{1}$-norm. We have that
\begin{equation*}
\begin{split}
\left\Vert \Delta _{j}\int_{0}^{t}T_{\Omega }(t-\tau )\nabla f(\tau )\ d\tau
\right\Vert _{L^{2}}& =\left\Vert \int_{0}^{t}T_{\Omega }(t-\tau )\nabla
\Delta _{j}f(\tau )\ d\tau \right\Vert _{L^{2}} \\
& \leq C\left\Vert \int_{0}^{t}e^{-(t-\tau )|\xi |^{2}}|\xi ||\widehat{\phi }%
_{j}(\xi )\widehat{f}(\tau )|\ d\tau \right\Vert _{L^{2}} \\
& \leq C\left\Vert \Vert e^{-(t-\tau )|\xi |^{2}}\Vert _{L_{\tau
}^{2}(0,t)}|\xi |\Vert \widehat{\phi }_{j}(\xi )\widehat{f}(\tau )\Vert
_{L_{\tau }^{2}(0,t)}\right\Vert _{L^{2}} \\
& \leq C\Vert \Delta _{j}f\Vert _{L^{2}(0,\infty ;L^{2})}.
\end{split}%
\end{equation*}%
Multiplying by $2^{sj}$, applying $l^{q}(\mathbb{Z})$-norm and using
inequality (\ref{embedding}), we arrive at
\begin{equation*}
\begin{split}
\left\Vert \int_{0}^{t}T_{\Omega }(t-\tau )\nabla f(\tau )\ d\tau
\right\Vert _{\dot{B}_{2,q}^{s}}& \leq C\left( \sum_{j\in \mathbb{Z}%
}2^{sjq}\Vert \Delta _{j}f\Vert _{L^{2}(0,\infty ;L^{2})}^{q}\right) ^{\frac{%
1}{q}} \\
& \leq C\Vert f\Vert _{L^{2}(0,\infty ;\dot{B}_{2,q}^{s})}
\end{split}%
\end{equation*}%
and then
\begin{equation}
\left\Vert \int_{0}^{t}T_{\Omega }(t-\tau )\nabla f(\tau )\ d\tau
\right\Vert _{X_{1}}\leq C\Vert f\Vert _{L^{2}(0,\infty ;\dot{B}_{2,q}^{s})}.
\label{3.21}
\end{equation}%
In order to estimate the $X_{2}$-norm, we use Lemma \ref{linearestimative1}
and Hardy-Littlewood-Sobolev inequality to obtain
\begin{equation}
\begin{split}
\left\Vert \int_{0}^{t}T_{\Omega }(t-\tau )\nabla f(\tau )\ d\tau
\right\Vert _{X_{2}}& \leq \left\Vert \int_{0}^{t}\Vert T_{\Omega }(t-\tau
)\nabla f(\tau )\Vert _{\dot{B}_{3,q}^{s}}\ d\tau \right\Vert
_{L^{4}(0,\infty )} \\
& \leq C\left\Vert \int_{0}^{t}(t-\tau )^{-\frac{1}{2}-\frac{3}{2}\left(
\frac{1}{2}-\frac{1}{3}\right) }\Vert f(\tau )\Vert _{\dot{B}_{2,q}^{s}}\
d\tau \right\Vert _{L^{4}(0,\infty )} \\
& \leq C\Vert f\Vert _{L^{2}(0,\infty ;\dot{B}_{2,q}^{s})}.
\end{split}
\label{3.22}
\end{equation}%
Putting together (\ref{3.21}) and (\ref{3.22}), we arrive at (\ref%
{aux-lemma-10}). \fin

\section{Global existence}

In this section we state and prove our results about existence and
uniqueness of global solutions to (\ref{NSC}). Basically, we have two cases $%
1/2<s<3/4$ and $s=1/2$. We start with the former.

\begin{theorem}
\label{theorem1}

\begin{enumerate}
\item[$(i)$] For $1\leq q<\infty $, consider $s,p$ and $\theta $ satisfying
\begin{gather*}
\frac{1}{2}<s<\frac{3}{4},\ \ \ \frac{1}{3}+\frac{s}{9}<\frac{1}{p}<\frac{2}{%
3}-\frac{s}{3}, \\
\frac{s}{2}-\frac{1}{2p}<\frac{1}{\theta }<\frac{5}{8}-\frac{3}{2p}+\frac{s}{%
4},\ \ \ \frac{3}{4}-\frac{3}{2p}\leq \frac{1}{\theta }<\min \left\{ 1-\frac{%
2}{p},\frac{1}{q}\right\} .
\end{gather*}%
Let $\Omega \in \mathbb{R}\setminus \{0\}$ and $u_{0}\in \dot{B}_{2,q}^{s}(%
\mathbb{R}^{3})$ with $\nabla \cdot u_{0}=0.$ There is a constant $%
C=C(s,p,\theta )>0$ such that if $\Vert u_{0}\Vert _{\dot{B}_{2,q}^{s}}\leq
C|\Omega |^{\frac{s}{2}-\frac{1}{4}}$, then there exists a unique global
solution $u\in C([0,\infty );\dot{B}_{2,q}^{s}(\mathbb{R}^{3}))$ to (\ref%
{NSC}).

\item[$(ii)$] For $q=\infty $, consider $s,p$ and $\theta $ satisfying
\begin{gather*}
\frac{1}{2}<s<\frac{3}{4},\ \ \ \frac{1}{3}+\frac{s}{9}<\frac{1}{p}<\frac{2}{%
3}-\frac{s}{3}, \\
\frac{s}{2}-\frac{1}{2p}<\frac{1}{\theta }<\frac{5}{8}-\frac{3}{2p}+\frac{s}{%
4}.
\end{gather*}%
Let $\Omega _{0}>0$ and $u_{0}\in \mathcal{I}$ with $\nabla \cdot u_{0}=0$,
where
\begin{equation}
\mathcal{I}:=\left\{ f\in \mathcal{S}^{\prime }(\mathbb{R}^{3})\colon \Vert
f\Vert _{\mathcal{I}}:=\sup_{|\Omega |\geq \Omega _{0}}|\Omega |^{\frac{1}{%
\theta }-\frac{3}{4}+\frac{3}{2p}}\Vert T_{\Omega }(t)f\Vert _{L^{\theta
}(0,\infty ;\dot{B}_{p,\infty }^{s})}<\infty \right\} .  \label{aux-space-77}
\end{equation}%
There is a constant $C=C(s,p,\theta )>0$ such that if $\Vert u_{0}\Vert _{%
\mathcal{I}}\leq C|\Omega |^{\frac{s}{2}-\frac{1}{4}}$ for $|\Omega |\geq
\Omega _{0}$, then the system (\ref{NSC}) has a unique global solution $u\in
L^{\theta }(0,\infty ;\dot{B}_{p,\infty }^{s}(\mathbb{R}^{3}))$. Moreover,
if in addition $u_{0}\in \dot{B}_{2,\infty }^{s}(\mathbb{R}^{3})$ then $u\in
C_{\omega }([0,\infty );\dot{B}_{2,\infty }^{s}(\mathbb{R}^{3}))$ where $%
C_{\omega }$ stands to time weakly continuous functions.
\end{enumerate}
\end{theorem}

\begin{remark}
Notice that the space $\mathcal{I}$ depends on the parameters $\Omega
_{0},\theta ,p$ and $s$, but for simplicity we have omitted them in the
notation.
\end{remark}

\bigskip

\textbf{Proof of Theorem \ref{theorem1}. }

Part $(i)$: By Lemma \ref{linearestimative1}, it follows that
\begin{equation}
\Vert T_{\Omega }(t)u_{0}\Vert _{L^{\theta }(0,\infty ;\dot{B}%
_{p,q}^{s})}\leq C_{0}|\Omega |^{-\frac{1}{\theta }+\frac{3}{4}-\frac{3}{2p}%
}\Vert u_{0}\Vert _{\dot{B}_{2,q}^{s}}.  \label{aux-linear-11}
\end{equation}%
Now, we define the operator $\Gamma $ and the set $Z$ by

\begin{equation}
\Gamma (u)(t)=T_{\Omega }(t)u_{0}-\mathfrak{B}(u,u)(t)  \label{operatorB}
\end{equation}%
and%
\begin{equation*}
Z=\left\{ u\in L^{\theta }(0,\infty ;\dot{B}_{p,q}^{s}(\mathbb{R}%
^{3})):\Vert u\Vert _{L^{\theta }(0,\infty ;\dot{B}_{p,q}^{s})}\leq
2C_{0}|\Omega |^{-\frac{1}{\theta }+\frac{3}{4}-\frac{3}{2p}}\Vert
u_{0}\Vert _{\dot{B}_{2,q}^{s}},\ \nabla \cdot u=0\right\} .
\end{equation*}%
Taking $\frac{1}{r}=\frac{2}{p}-\frac{s}{3},$ we can employ Lemma \ref%
{nonlinear1} and (\ref{remark1}) to estimate $\Gamma (\cdot )$ as follows
\begin{equation}
\begin{split}
\Vert \Gamma (u)-\Gamma (v)& \Vert _{L^{\theta }(0,\infty ;\dot{B}%
_{p,q}^{s})}=\left\Vert \int_{0}^{t}T_{\Omega }(t-\tau )\mathbb{P}\nabla
\cdot (u\otimes (u-v)(\tau )+(u-v)\otimes v(\tau ))\ \tau \right\Vert
_{L^{\theta }(0,\infty ;\dot{B}_{p,q}^{s})} \\
& \leq C|\Omega |^{\frac{1}{\theta }-\frac{1}{2}+\frac{3}{2}\left( \frac{1}{r%
}-\frac{1}{p}\right) }\Vert u\otimes (u-v)+(u-v)\otimes v\Vert _{L^{\frac{%
\theta }{2}}(0,\infty ;\dot{B}_{r,q}^{s})} \\
& \leq C|\Omega |^{\frac{1}{\theta }-\frac{1}{2}+\frac{3}{2}\left( \frac{1}{r%
}-\frac{1}{p}\right) }\left( \Vert u\Vert _{L^{\theta }(0,\infty ;\dot{B}%
_{p,q}^{s})}+\Vert v\Vert _{L^{\theta }(0,\infty ;\dot{B}_{p,q}^{s})}\right)
\Vert u-v\Vert _{L^{\theta }(0,\infty ;\dot{B}_{p,q}^{s})} \\
& \leq C|\Omega |^{\frac{1}{\theta }-\frac{1}{2}+\frac{3}{2}\left( \frac{1}{r%
}-\frac{1}{p}\right) }4C_{0}|\Omega |^{-\frac{1}{\theta }+\frac{3}{4}-\frac{3%
}{2p}}\Vert u_{0}\Vert _{\dot{B}_{2,q}^{s}}\Vert u-v\Vert _{L^{\theta
}(0,\infty ;\dot{B}_{p,q}^{s})} \\
& =C_{2}|\Omega |^{\frac{1}{\theta }-\frac{1}{2}+\frac{3}{2}\left( \frac{1}{r%
}-\frac{1}{p}\right) -\frac{1}{\theta }+\frac{3}{4}-\frac{3}{2p}}\Vert
u_{0}\Vert _{\dot{B}_{2,q}^{s}}\Vert u-v\Vert _{L^{\theta }(0,\infty ;\dot{B}%
_{p,q}^{s})} \\
& =C_{2}|\Omega |^{\frac{1}{4}-\frac{s}{2}}\Vert u_{0}\Vert _{\dot{B}%
_{2,q}^{s}}\Vert u-v\Vert _{L^{\theta }(0,\infty ;\dot{B}_{p,q}^{s})}, \\
&
\end{split}
\label{contraction}
\end{equation}%
for all $u,v\in Z$, where $C_{2}=C_{2}(s,p,\theta ).$ Moreover, using (\ref%
{aux-linear-11}) and (\ref{contraction}) with $v=0$, we obtain
\begin{equation}
\begin{split}
\Vert \Gamma (u)\Vert _{L^{\theta }(0,\infty ;\dot{B}_{p,q}^{s})}& \leq
\Vert T_{\Omega }(t)u_{0}\Vert _{L^{\theta }(0,\infty ;\dot{B}%
_{p,q}^{s})}+\Vert \Gamma (u)-\Gamma (0)\Vert _{L^{\theta }(0,\infty ;\dot{B}%
_{p,q}^{s})} \\
& \leq C_{0}|\Omega |^{-\frac{1}{\theta }+\frac{3}{4}-\frac{3}{2p}}\Vert
u_{0}\Vert _{\dot{B}_{2,q}^{s}}+C_{2}|\Omega |^{\frac{1}{4}-\frac{s}{2}%
}\Vert u_{0}\Vert _{\dot{B}_{2,q}^{s}}\Vert u\Vert _{L^{\theta }(0,\infty ;%
\dot{B}_{p,q}^{s})} \\
& \leq C_{0}|\Omega |^{-\frac{1}{\theta }+\frac{3}{4}-\frac{3}{2p}}\Vert
u_{0}\Vert _{\dot{B}_{2,q}^{s}}+C_{2}|\Omega |^{\frac{1}{4}-\frac{s}{2}%
}\Vert u_{0}\Vert _{\dot{B}_{2,q}^{s}}2C_{0}|\Omega |^{-\frac{1}{\theta }+%
\frac{3}{4}-\frac{3}{2p}}\Vert u_{0}\Vert _{\dot{B}_{2,q}^{s}} \\
& =C_{0}\Vert u_{0}\Vert _{\dot{B}_{2,q}^{s}}|\Omega |^{-\frac{1}{\theta }+%
\frac{3}{4}-\frac{3}{2p}}\left( 1+2C_{2}|\Omega |^{\frac{1}{4}-\frac{s}{2}%
}\Vert u_{0}\Vert _{\dot{B}_{2,q}^{s}}\right)
\end{split}
\label{differenceofoperator}
\end{equation}%
for all $u\in Z$. Thus, for $\Omega $ and $u_{0}$ satisfying
\begin{equation*}
C_{2}|\Omega |^{\frac{1}{4}-\frac{s}{2}}\Vert u_{0}\Vert _{\dot{B}%
_{2,q}^{s}}\leq \frac{1}{2},
\end{equation*}%
we get
\begin{equation*}
\Vert \Gamma (u)\Vert _{L^{\theta }(0,\infty ;\dot{B}_{p,q}^{s})}\leq
2C_{0}|\Omega |^{-\frac{1}{\theta }+\frac{3}{4}-\frac{3}{2p}}\Vert
u_{0}\Vert _{\dot{B}_{2,q}^{s}}\text{ and }\Vert \Gamma (u)-\Gamma (v)\Vert
_{L^{\theta }(0,\infty ;\dot{B}_{p,q}^{s})}\leq \frac{1}{2}\Vert u-v\Vert
_{L^{\theta }(0,\infty ;\dot{B}_{p,q}^{s})}.
\end{equation*}%
Then, Banach fixed point theorem implies that there exists a unique mild
solution $u\in Z$ to (\ref{NSC}), i.e.,
\begin{equation*}
u(t)=T_{\Omega }(t)u_{0}-\mathfrak{B}(u,u)(t).
\end{equation*}%
It remains to prove that $u\in C([0,\infty );\dot{B}_{2,q}^{s}(\mathbb{R}%
^{3}))$. Basically, we need to estimate the $\dot{B}_{2,q}^{s}$-norm of the
linear and nonlinear parts in (\ref{operatorB}). For the linear one, we use
Lemma \ref{linearestimative1} to get
\begin{equation}
\Vert T_{\Omega }(t)u_{0}\Vert _{\dot{B}_{2,q}^{s}}\leq C_{0}\Vert
u_{0}\Vert _{\dot{B}_{2,q}^{s}}.  \label{4.9}
\end{equation}%
For the nonlinear part, taking $\frac{1}{r}=\frac{2}{p}-\frac{s}{3},$ we use
Lemma \ref{linearestimative1}, (\ref{remark1}) and H\"{o}lder inequality to
obtain

\begin{equation}
\begin{split}
\left\Vert \int_{0}^{t}T_{\Omega }(t-\tau )\mathbb{P}\nabla \cdot (u\otimes
u)(\tau )\ d\tau \right\Vert _{\dot{B}_{2,q}^{s}}& \leq
C\int_{0}^{t}\left\Vert T_{\Omega }(t-\tau )\mathbb{P}\nabla \cdot (u\otimes
u)(\tau )\right\Vert _{\dot{B}_{2,q}^{s}}\ d\tau \\
& \leq C\int_{0}^{t}(t-\tau )^{-\frac{1}{2}-\frac{3}{2}\left( \frac{1}{r}-%
\frac{1}{2}\right) }\left\Vert (u\otimes u)(\tau )\right\Vert _{\dot{B}%
_{r,q}^{s}}\ d\tau \\
& \leq C\int_{0}^{t}(t-\tau )^{-\frac{1}{2}-\frac{3}{2}\left( \frac{1}{r}-%
\frac{1}{2}\right) }\left\Vert u(\tau )\right\Vert _{\dot{B}_{p,q}^{s}}^{2}\
d\tau \\
& \leq C\left\Vert (t-\cdot )^{-\frac{1}{2}-\frac{3}{2r}+\frac{3}{4}%
}\right\Vert _{L^{\frac{\theta }{\theta -2}}(0<\tau <t)}\left\Vert \Vert
u(\tau )\Vert _{\dot{B}_{p,q}^{s}}^{2}\right\Vert _{L^{\frac{\theta }{2}%
}(0,\infty )} \\
& \leq Ct^{\frac{\theta -2}{\theta }\left( 1+\frac{\theta }{\theta -2}\left(
-\frac{1}{2}-\frac{3}{2r}+\frac{3}{4}\right) \right) }\Vert u\Vert
_{L^{\theta }(0,\infty ;\dot{B}_{p,q}^{s})}^{2},
\end{split}
\label{4.10}
\end{equation}%
where we need $\frac{1}{\theta }<\frac{5}{8}-\frac{3}{2p}+\frac{s}{4}$ in
order to ensure integrability at $\tau =t$. From (\ref{4.9}) and (\ref{4.10}%
), it follows that $u(t)\in \dot{B}_{2,q}^{s}(\mathbb{R}^{3})$ for $t>0$,
and then we have that $u\in C([0,\infty );\dot{B}_{2,q}^{s}(\mathbb{R}^{3}))$%
, as desired. 

Part $(ii)$: In view of (\ref{aux-space-77}), we have that
\begin{equation}
\Vert T_{\Omega }(t)u_{0}\Vert _{L^{\theta }(0,\infty ;\dot{B}_{p,\infty
}^{s})}\leq |\Omega |^{-\frac{1}{\theta }+\frac{3}{4}-\frac{3}{2p}}\Vert
u_{0}\Vert _{\mathcal{I}},\text{ for all }\left\vert \Omega \right\vert \geq
\Omega _{0}\text{.}  \label{aux-linear-11_q=infty}
\end{equation}%
Now, for $\left\vert \Omega \right\vert \geq \Omega _{0}$ consider

\begin{equation}
\Gamma (u)(t)=T_{\Omega }(t)u_{0}-\mathfrak{B}(u,u)(t)
\label{operatorB_q=infty}
\end{equation}%
and%
\begin{equation*}
Z=\left\{ u\in L^{\theta }(0,\infty ;\dot{B}_{p,\infty }^{s}(\mathbb{R}%
^{3})):\Vert u\Vert _{L^{\theta }(0,\infty ;\dot{B}_{p,\infty }^{s})}\leq
2|\Omega |^{-\frac{1}{\theta }+\frac{3}{4}-\frac{3}{2p}}\Vert u_{0}\Vert _{%
\mathcal{I}},\ \nabla \cdot u=0\right\} .
\end{equation*}%
Taking $\frac{1}{r}=\frac{2}{p}-\frac{s}{3},$ and proceeding similarly to
Part $(i)$, we obtain a constant $\tilde{C_{2}}=\tilde{C_{2}}(s,p,\theta )$
such that
\begin{equation}
\begin{split}
\Vert \Gamma (u)-\Gamma (v)\Vert _{L^{\theta }(0,\infty ;\dot{B}_{p,\infty
}^{s})}& \leq \tilde{C_{2}}|\Omega |^{\frac{1}{4}-\frac{s}{2}}\Vert
u_{0}\Vert _{\mathcal{I}}\Vert u-v\Vert _{L^{\theta }(0,\infty ;\dot{B}%
_{p,\infty }^{s})} \\
\Vert \Gamma (u)\Vert _{L^{\theta }(0,\infty ;\dot{B}_{p,\infty }^{s})}&
\leq \Vert u_{0}\Vert _{\mathcal{I}}|\Omega |^{-\frac{1}{\theta }+\frac{3}{4}%
-\frac{3}{2p}}\left( 1+2\tilde{C_{2}}|\Omega |^{\frac{1}{4}-\frac{s}{2}%
}\Vert u_{0}\Vert _{\mathcal{I}}\right) ,
\end{split}
\label{contraction-2222}
\end{equation}%
for all $u,v\in Z$. Thus, for $\Omega $ and $u_{0}$ satisfying
\begin{equation*}
\left\vert \Omega \right\vert \geq \Omega _{0}\text{ and }\tilde{C_{2}}%
|\Omega |^{\frac{1}{4}-\frac{s}{2}}\Vert u_{0}\Vert _{\mathcal{I}}\leq \frac{%
1}{2} ,
\end{equation*}%
we get
\begin{equation*}
\Vert \Gamma (u)\Vert _{L^{\theta }(0,\infty ;\dot{B}_{p,\infty }^{s})}\leq
2|\Omega |^{-\frac{1}{\theta }+\frac{3}{4}-\frac{3}{2p}}\Vert u_{0}\Vert _{%
\mathcal{I}}\text{ and }\Vert \Gamma (u)-\Gamma (v)\Vert _{L^{\theta
}(0,\infty ;\dot{B}_{p,\infty }^{s})}\leq \frac{1}{2}\Vert u-v\Vert
_{L^{\theta }(0,\infty ;\dot{B}_{p,\infty }^{s})}.
\end{equation*}%
Again, we can apply the Banach fixed point theorem in order to obtain a
unique mild solution $u\in Z$ to (\ref{NSC}). Assume now that $u_{0}\in \dot{%
B}_{2,\infty }^{s}(\mathbb{R}^{3})$. Since (\ref{4.9}) and (\ref{4.10}) hold
true for $q=\infty $, it follows that $u\in C_{\omega }([0,\infty );\dot{B}%
_{2,\infty }^{s}(\mathbb{R}^{3}))$. \fin

\bigskip

Before proceeding, for $\Omega _{0}>0$ and $1\leq q\leq \infty $ we define
the space
\begin{equation}
\mathcal{F}:=\left\{ f\in \mathcal{S}^{\prime }(\mathbb{R}^{3})\colon \Vert
f\Vert _{\mathcal{F}}:=\sup_{|\Omega |\geq \Omega _{0}}\Vert T_{\Omega
}(t)f\Vert _{L^{4}(0,\infty ;\dot{B}_{3,q}^{1/2})}<\infty \right\} ,
\label{aux-space-771}
\end{equation}%
where, for simplicity, we have omitted the dependence on $\Omega _{0}$ and $%
q $ in the notation $\mathcal{F}$. We also define
\begin{equation}
\mathcal{F}_{0}:=\left\{ f\in \mathcal{F}\colon \limsup_{|\Omega
|\rightarrow \infty }\Vert T_{\Omega }(t)f\Vert _{L^{4}(0,\infty ;\dot{B}%
_{3,q}^{1/2})}=0\right\} .  \label{aux-space-772}
\end{equation}%
Both spaces $\mathcal{F}$ and $\mathcal{F}_{0}$ are endowed with the norm $%
\Vert \cdot \Vert _{\mathcal{F}}.$ The next theorem deals with the critical
case $s=1/2$.

\begin{theorem}
\label{theorem3}Let $2\leq q\leq \infty $ and $u_{0}\in D$ with $\nabla
\cdot u_{0}=0$ where $D$ is a precompact set in $\mathcal{F}_{0}$. Then,
there exist $\widetilde{\Omega }=\widetilde{\Omega }(D)>0$ and a unique
global solution $u$ to (\ref{NSC}) in $L^{4}(0,\infty ;\dot{B}_{3,q}^{1/2}(%
\mathbb{R}^{3}))$ provided that $|\Omega |\geq \widetilde{\Omega }$.
Moreover, if in addition $u_{0}\in \dot{B}_{2,q}^{1/2}(\mathbb{R}^{3})$ with
$q\neq \infty $, then $u\in C([0,\infty );\dot{B}_{2,q}^{1/2}(\mathbb{R}%
^{3}))$. In the case $q=\infty ,$ we obtain $C_{\omega }([0,\infty );\dot{B}%
_{2,\infty }^{1/2}(\mathbb{R}^{3})).$
\end{theorem}

\textbf{Proof. } Let $\delta $ be a positive number that will be chosen
later. Given that $D$ is a precompact set in $\mathcal{F}_{0}$, there exist $%
L=L(\delta ,D)\in \mathbb{N}$ and $\{g_{k}\}\subset \mathcal{F}_{0}$ such
that
\begin{equation*}
D\subset \bigcup_{k=1}^{L}B(g_{k},\delta ),
\end{equation*}%
where $B(g_{k},\delta )$ denotes the ball in $\mathcal{F}_{0}$ with center $%
g_{k}$ and radius $\delta $. On the other hand, using the definition (\ref%
{aux-space-772}), there exists $\tilde{\Omega}=\tilde{\Omega}(\delta ,D)\geq
\Omega _{0}>0$ such that
\begin{equation*}
\sup_{k=1,2,\ldots ,L}\Vert T_{\Omega }(t)g_{k}\Vert _{L^{4}(0,\infty ;\dot{B%
}_{3,q}^{\frac{1}{2}})}\leq \delta
\end{equation*}%
provided that $|\Omega |\geq \tilde{\Omega}$. Now, given $g\in D$ there
exists $k\in \{1,2,\ldots ,L\}$ such that $g\in B(g_{k},\delta )$.
Therefore, for $|\Omega |\geq \tilde{\Omega}$ we can estimate
\begin{equation*}
\begin{split}
\Vert T_{\Omega }(t)g\Vert _{L^{4}(0,\infty ;\dot{B}_{3,q}^{\frac{1}{2}})}&
\leq \Vert T_{\Omega }(t)(g_{k}-g)\Vert _{L^{4}(0,\infty ;\dot{B}_{3,q}^{%
\frac{1}{2}})}+\Vert T_{\Omega }(t)g_{k}\Vert _{L^{4}(0,\infty ;\dot{B}%
_{3,q}^{\frac{1}{2}})} \\
& \leq C\Vert g_{k}-g\Vert _{\mathcal{F}}+\delta \\
& \leq (C+1)\delta .
\end{split}%
\end{equation*}%
Thus, there exists $C_{1}>0$ such that
\begin{equation}
\sup_{g\in D}\Vert T_{\Omega }(t)g\Vert _{L^{4}(0,\infty ;\dot{B}_{3,q}^{%
\frac{1}{2}})}\leq C_{1}\delta ,\text{ for all }|\Omega |\geq \tilde{\Omega}.
\label{aux-3.5}
\end{equation}%
Now, we consider the complete metric space $Z$ defined by
\begin{equation}
Z=\left\{ u\in \ L^{4}(0,\infty ;\dot{B}_{3,q}^{\frac{1}{2}}):\Vert u\Vert
_{L^{4}(0,\infty ;\dot{B}_{3,q}^{\frac{1}{2}})}\leq 2C_{1}\delta ,\nabla
\cdot u=0\right\} ,\   \label{space-1003}
\end{equation}%
endowed with the metric $d(u,v)=\Vert u-v\Vert _{L^{4}(0,\infty ;\dot{B}%
_{3,q}^{\frac{1}{2}})}.$ Also, we consider the operator $\Gamma $ defined in
the proof of Theorem \ref{theorem1}. For $u,v\in Z$, using Lemma \ref%
{nonlinearestimativecritical}, (\ref{remark1}) and H\"{o}lder inequality, we
can estimate
\begin{equation}
\begin{split}
\Vert \Gamma (u)-\Gamma (v)\Vert _{L^{4}(0,\infty ;\dot{B}_{3,q}^{\frac{1}{2}%
})}& =\left\Vert \int_{0}^{t}T_{\Omega }(t-\tau )\mathbb{P}\nabla \cdot
(u\otimes (u-v)+(u-v)\otimes v)(\tau )\ d\tau \right\Vert _{L^{4}(0,\infty ;%
\dot{B}_{3,q}^{\frac{1}{2}})} \\
& \leq C\Vert u\otimes (u-v)+(u-v)\otimes v\Vert _{L^{2}(0,\infty ;\dot{B}%
_{2,q}^{\frac{1}{2}})} \\
& \leq C\left( \left\Vert \Vert u\Vert _{\dot{B}_{3,q}^{\frac{1}{2}}}\Vert
u-v\Vert _{\dot{B}_{3,q}^{\frac{1}{2}}}\right\Vert _{L^{2}(0,\infty
)}+\left\Vert \Vert v\Vert _{\dot{B}_{3,q}^{\frac{1}{2}}}\Vert u-v\Vert _{%
\dot{B}_{3,q}^{\frac{1}{2}}}\right\Vert _{L^{2}(0,\infty )}\right) \\
& \leq C_{2}\left( \Vert u\Vert _{L^{4}(0,\infty ;\dot{B}_{3,q}^{\frac{1}{2}%
})}+\Vert v\Vert _{L^{4}(0,\infty ;\dot{B}_{3,q}^{\frac{1}{2}})}\right)
\Vert u-v\Vert _{L^{4}(0,\infty ;\dot{B}_{3,q}^{\frac{1}{2}})}.
\end{split}
\label{3.6}
\end{equation}%
Taking $v=0$ in (\ref{3.6}), for $u\in Z$ it follows that
\begin{equation}
\begin{split}
\Vert \Gamma (u)\Vert _{L^{4}(0,\infty ;\dot{B}_{3,q}^{\frac{1}{2}})}& \leq
\Vert \Gamma (0)\Vert _{L^{4}(0,\infty ;\dot{B}_{3,q}^{\frac{1}{2}})}+\Vert
\Gamma (u)-\Gamma (0)\Vert _{L^{4}(0,\infty ;\dot{B}_{3,q}^{\frac{1}{2}})} \\
& \leq \Vert T_{\Omega }(t)u_{0}\Vert _{L^{4}(0,\infty ;\dot{B}_{3,q}^{\frac{%
1}{2}})}+C_{2}\Vert u\Vert _{L^{4}(0,\infty ;\dot{B}_{3,q}^{\frac{1}{2}%
})}^{2}.
\end{split}
\label{3.7}
\end{equation}%
Choosing $0<\delta <\frac{1}{8C_{1}C_{2}}$, estimates (\ref{aux-3.5}), (\ref%
{3.6}) and (\ref{3.7}) yield
\begin{equation*}
\begin{split}
& \Vert \Gamma (u)\Vert _{L^{4}(0,\infty ;\dot{B}_{3,q}^{\frac{1}{2}})}\leq
2C_{1}\delta ,\text{ for all }u\in Z, \\
& \Vert \Gamma (u)-\Gamma (v)\Vert _{L^{4}(0,\infty ;\dot{B}_{3,q}^{\frac{1}{%
2}})}\leq \frac{1}{2}\Vert u-v\Vert _{L^{4}(0,\infty ;\dot{B}_{3,q}^{\frac{1%
}{2}})},\text{ for all }u,v\in Z,
\end{split}%
\end{equation*}%
provided that $|\Omega |\geq \tilde{\Omega}$. Therefore, we can apply the
Banach fixed point theorem to obtain a unique global solution $u\in \
L^{4}(0,\infty ;\dot{B}_{3,q}^{\frac{1}{2}}).$

Moreover, using Lemma \ref{linearestimative1}, Lemma \ref%
{nonlinearestimativecritical} and $u\in L^{4}(0,\infty ;\dot{B}_{3,q}^{\frac{%
1}{2}})$, we have that

\begin{equation}
\Vert u(t)\Vert _{\dot{B}_{2,q}^{\frac{1}{2}}}=\Vert \Gamma (u)(t)\Vert _{%
\dot{B}_{2,q}^{\frac{1}{2}}}\leq C\Vert u_{0}\Vert _{\dot{B}_{2,q}^{\frac{1}{%
2}}}+C\Vert u\Vert _{L^{4}(0,\infty ;\dot{B}_{3,q}^{\frac{1}{2}%
})}^{2}<\infty ,  \label{aux-3.8}
\end{equation}%
for a.e. $t>0$. Since $\Vert u\Vert _{L^{4}(0,\infty ;\dot{B}_{3,q}^{\frac{1%
}{2}})}\leq 2C_{1}\delta $ and $\delta <\frac{1}{8C_{1}C_{2}}$, it follows
that
\begin{equation}
\Vert u(t)\Vert _{\dot{B}_{2,q}^{\frac{1}{2}}}\leq C(\Vert u_{0}\Vert _{\dot{%
B}_{2,q}^{\frac{1}{2}}}+1)<\infty ,\text{ for all }\left\vert \Omega
\right\vert \geq \tilde{\Omega},  \label{aux-3.9}
\end{equation}%
and so $u(t)\in \dot{B}_{2,q}^{\frac{1}{2}}(\mathbb{R}^{3})$ for a.e. $t>0$.
Using this and above estimates, standard arguments yield $u\in C([0,\infty );%
\dot{B}_{2,q}^{\frac{1}{2}}(\mathbb{R}^{3}))$ for $q\neq \infty $ and $u\in
C_{\omega }([0,\infty );\dot{B}_{2,q}^{\frac{1}{2}}(\mathbb{R}^{3}))$ for $%
q=\infty $.\fin

\begin{theorem}
\label{theorem2}Let $2\leq q\leq \infty $ and $u_{0}\in \mathcal{F}_{0}$
with $\nabla \cdot u_{0}=0$. Then, there exist $\tilde{\Omega}=\tilde{\Omega}%
(u_{0})$ and a unique global solution $u\in L^{4}(0,\infty ;\dot{B}_{3,q}^{%
\frac{1}{2}}(\mathbb{R}^{3}))$ to (\ref{NSC}) provided that $|\Omega |\geq
\tilde{\Omega}$.
\end{theorem}

\textbf{Proof. } It is sufficient to apply Theorem \ref{theorem3} to the set
$D=\{u_{0}\}$. \fin

\begin{cor}
\label{cor1} Let $2\leq q<4$ and $u_{0}\in D$ with $\nabla \cdot u_{0}=0$
where $D$ is a precompact set in $\dot{B}_{2,q}^{\frac{1}{2}}(\mathbb{R}%
^{3}) $. Then, there exist $\tilde{\Omega}(D)>0$ and a unique global
solution $u$ to (\ref{NSC}) in the class $C([0,\infty );\dot{B}_{2,q}^{\frac{%
1}{2}}(\mathbb{R}^{3}))\cap L^{4}(0,\infty ;\dot{B}_{3,q}^{\frac{1}{2}}(%
\mathbb{R}^{3}))$ provided that $|\Omega |\geq \tilde{\Omega}(D)$.
\end{cor}

\textbf{Proof. } In view of Lemma \ref{critical_linear_semigroup}, we have
that $\dot{B}_{2,q}^{\frac{1}{2}}\hookrightarrow \mathcal{F}_{0}$ for $1\leq
q<4$. Now the result follows by applying Theorem \ref{theorem3}.\fin

\section{Asymptotic behavior as $|\Omega |\rightarrow \infty $}

In this section we study the asymptotic behavior of the mild solutions as $%
|\Omega |\rightarrow \infty $. For convenience, we denote
\begin{equation*}
\alpha _{0}=-\frac{1}{\theta }+\frac{1}{2}-\frac{3}{2p}+\frac{s}{2}\ \ \text{
and }\ \ \beta _{0}=\frac{1}{\theta }-\frac{3}{4}+\frac{3}{2p}.
\end{equation*}

First, we consider the case $1/2<s<3/4$.

\begin{theorem}
\label{properties_of_u_and_v_Omega}

\begin{enumerate}
\item[$(i)$] Let $0\leq \epsilon <\frac{1}{12}$ and $1\leq q< \infty$, and
suppose that $s,p$ and $\theta $ satisfy%
\begin{gather*}
\frac{1}{2}+3\epsilon <s<\frac{3}{4},\ \ \ \ \frac{1}{3}+\frac{s}{9}<\frac{1%
}{p}<\frac{2}{3}-\frac{s}{3}, \\
\frac{s}{2}-\frac{1}{2p}<\frac{1}{\theta }<\frac{5}{8}-\frac{3}{2p}+\frac{s}{%
4}-\frac{\epsilon }{4},\ \ \ \ \frac{3}{4}-\frac{3}{2p}\leq \frac{1}{\theta }%
<\min\left\{1-\frac{2}{p},\frac{1}{q}\right\}.
\end{gather*}%
Let $u$ and $v$ be solutions of (\ref{NSC}) with initial data $u_{0}$ and $%
v_{0}$ in $\dot{B}_{2,q}^{s}(\mathbb{R}^{3})$, respectively. Then, for $%
\alpha <2\beta _{0}$
\begin{equation}
\lim_{|\Omega |\rightarrow \infty }|\Omega |^{\alpha }\Vert u(t)-v(t)\Vert _{%
\dot{B}_{2,q}^{s+\epsilon }}=0\ \ \ \text{if and only if}\ \ \ \lim_{|\Omega
|\rightarrow \infty }|\Omega |^{\alpha }\Vert T_{\Omega
}(t)(u_{0}-v_{0})\Vert _{\dot{B}_{2,q}^{s+\epsilon }}=0,\text{ for each
fixed }t>0.  \label{asymp-200}
\end{equation}

\item[$(ii)$] Let $0\leq \epsilon <\frac{1}{6}$ and $1\leq q< \infty$.
Assume that $s$, $p$ and $\theta $ satisfy
\begin{gather*}
\frac{1}{2}+\frac{3\epsilon }{2}<s<\frac{3}{4},\ \ \ \ \frac{1}{3}+\frac{s}{9%
}<\frac{1}{p}<\frac{2}{3}-\frac{s}{3}, \\
\frac{s}{2}-\frac{1}{2p}<\frac{1}{\theta }<\frac{5}{8}-\frac{3}{2p}+\frac{s}{%
4},\ \ \ \ \frac{3}{4}-\frac{3}{2p}\leq \frac{1}{\theta }<\min\left\{1-\frac{%
2}{p},\frac{1}{q}\right\}.
\end{gather*}%
Let $\alpha <\alpha _{0}+2\beta _{0}-\frac{\epsilon }{2}$ and assume that $u$
and $v$ are solutions of (\ref{NSC}) with initial data $u_{0}$ and $v_{0}$
in $\dot{B}_{2,q}^{s}(\mathbb{R}^{3})$, respectively. Then, for each fixed $%
t>0$,
\begin{equation}
\lim_{|\Omega |\rightarrow \infty }|\Omega |^{\alpha }\Vert u-v\Vert _{{%
L^{\theta }(0,\infty ;\dot{B}_{p,q}^{s+\epsilon })}}=0\ \ \ \text{if and
only if}\ \ \ \lim_{|\Omega |\rightarrow \infty }|\Omega |^{\alpha }\Vert
T_{\Omega }(t)(u_{0}-v_{0})\Vert _{L^{\theta }(0,\infty ;\dot{B}%
_{p,q}^{s+\epsilon })}=0.  \label{asymp-201}
\end{equation}
\end{enumerate}
\end{theorem}

\textbf{Proof. } First we write

\begin{equation}
u-v=T_{\Omega }(t)(u_{0}-v_{0})+\mathfrak{B}(u,u)(t)-\mathfrak{B}(v,v)(t).
\label{aux-1}
\end{equation}

Considering $\frac{1}{r}=\frac{2}{p}-\frac{s}{3}$, we estimate the $\dot{B}%
_{2,q}^{s+\epsilon }$-norm of the nonlinear term in (\ref{aux-1}) as follows
\begin{equation*}
\begin{split}
\Vert \mathfrak{B}(u,u)(t)-\mathfrak{B}(v,v)(t)\Vert _{\dot{B}%
_{2,q}^{s+\epsilon }}& \leq C\int_{0}^{t}(t-\tau )^{-\frac{1}{2}-\frac{3}{2}%
\left( \frac{1}{r}-\frac{1}{2}\right) }\Vert e^{\frac{1}{2}(t-\tau )\Delta
}(u\otimes (u-v)+(u-v)\otimes v)(\tau )\Vert _{\dot{B}_{r,q}^{s+\epsilon }}\
d\tau \\
& \leq C\int_{0}^{t}(t-\tau )^{-\frac{1}{2}-\frac{3}{2}\left( \frac{1}{r}-%
\frac{1}{2}\right) -\frac{\epsilon }{2}}(\Vert u(\tau )\Vert _{\dot{B}%
_{p,q}^{s}}+\Vert u(\tau )\Vert _{\dot{B}_{p,q}^{s}})\Vert (u-v)(\tau )\Vert
_{\dot{B}_{p,q}^{s}}\ d\tau \\
& \leq Ct^{\frac{1}{2}-\frac{2}{\theta }-\frac{3}{2}\left( \frac{1}{r}-\frac{%
1}{2}\right) -\frac{\epsilon }{2}}(\Vert u\Vert _{L^{\theta }(0,\infty ;\dot{%
B}_{p,q}^{s})}+\Vert v\Vert _{L^{\theta }(0,\infty ;\dot{B}%
_{p,q}^{s})})\Vert u-v\Vert _{L^{\theta }(0,\infty ;\dot{B}_{p,q}^{s})},
\end{split}%
\end{equation*}%
where we have the integrability at $\tau =t$ due to the condition
\begin{equation*}
\frac{1}{\theta }<\frac{5}{8}-\frac{3}{2p}+\frac{s}{4}-\frac{\epsilon }{4}%
\Longrightarrow \frac{1}{2}-\frac{2}{\theta }-\frac{3}{2}\left( \frac{1}{r}-%
\frac{1}{2}\right) -\frac{\epsilon }{2}>0.
\end{equation*}%
Thus,
\begin{equation*}
|\Omega |^{\alpha }\Vert \mathfrak{B}(u,u)(t)-\mathfrak{B}(v,v)(t)\Vert _{%
\dot{B}_{2,q}^{s+\epsilon }}\leq Ct^{\frac{1}{2}-\frac{2}{\theta }-\frac{3}{2%
}\left( \frac{1}{r}-\frac{1}{2}\right) -\frac{\epsilon }{2}}|\Omega
|^{\alpha -2\beta _{0}},\text{ for all }|\Omega |\geq \Omega _{0}.
\end{equation*}%
Since $\alpha <2\beta _{0}$, it follows that
\begin{equation}
\lim_{|\Omega |\rightarrow \infty }|\Omega |^{\alpha }\Vert \mathfrak{B}%
(u,u)(t)-\mathfrak{B}(v,v)(t)\Vert _{\dot{B}_{2,q}^{s+\epsilon }}=0,\text{
for each }t>0.  \label{aux-101}
\end{equation}%
In view of (\ref{aux-1}) and (\ref{aux-101}), we obtain the desired property.

For item $(ii)$, we proceed similarly as in the proof of Lemma \ref%
{nonlinear1} by taking $f=u\otimes (u-v)+(u-v)\otimes v$ in the nonlinear
term of (\ref{aux-1}). Since $\frac{1}{\theta }<\frac{1}{2}-\frac{3}{2}%
\left( \frac{1}{r}-\frac{1}{p}\right) -\frac{\epsilon }{2},$ we can estimate
\begin{equation*}
\Vert \mathfrak{B}(u,u)-\mathfrak{B}(v,v)\Vert _{L^{\theta }(0,\infty ;\dot{B%
}_{p,q}^{s+\epsilon })}\leq C|\Omega |^{-\alpha _{0}+\frac{\epsilon }{2}%
}(\Vert u\Vert _{L^{\theta }(0,\infty ;\dot{B}_{p,q}^{s})}+\Vert v\Vert
_{L^{\theta }(0,\infty ;\dot{B}_{p,q}^{s})})\Vert u-v\Vert _{L^{\theta
}(0,\infty ;\dot{B}_{p,q}^{s})}
\end{equation*}%
and then
\begin{equation}
|\Omega |^{\alpha }\Vert \mathfrak{B}(u,u)-\mathfrak{B}(v,v)\Vert
_{L^{\theta }(0,\infty ;\dot{B}_{p,q}^{s+\epsilon })}\leq C|\Omega |^{\alpha
-\alpha _{0}+\frac{\epsilon }{2}-2\beta _{0}},\text{ for all }|\Omega |\geq
\Omega _{0}.  \label{aux-201}
\end{equation}%
Finally, we obtain (\ref{asymp-201}) by letting $|\Omega |\rightarrow \infty
$ and using (\ref{aux-201}) and (\ref{aux-1}).\fin

\begin{remark}
Let $1\leq q<\infty $, and consider $s,\gamma _{2},p$ and $\theta $ such
that
\begin{gather}
\frac{1}{2}<s<\frac{3}{4},\ \ \ 0<\gamma _{2}<\frac{1}{2}\left( 1-\frac{1}{%
\epsilon }\right) ,\ \ \ \frac{1}{2\gamma _{2}}\left( \frac{1}{8}-\frac{s}{4}%
+\gamma _{2}\right) \leq \frac{1}{p}<\frac{2}{3}-\frac{s}{3},
\label{new_assumptions1} \\
\frac{s}{2}-\frac{1}{2p}<\frac{1}{\theta }<\frac{5}{8}-\frac{3}{2p}+\frac{s}{%
4},\ \ \ \frac{3}{4}-\frac{3}{2p}\leq \frac{1}{\theta }<1-\frac{2}{p}.
\label{new_assumptions2}
\end{gather}%
Since $\frac{1}{\theta }<\frac{1}{2}-\frac{3}{2}\left( \frac{1}{r}-\frac{1}{p%
}\right) -\gamma _{2}\left( 1-\frac{2}{p}\right) $ and
\begin{equation*}
\left( \frac{\log {(e+|\Omega |t)}}{1+|\Omega |t}\right) ^{\frac{1}{2}\left(
1-\frac{2}{p}\right) }=\frac{\log {(e+|\Omega |t)^{\frac{1}{2}\left( 1-\frac{%
2}{p}\right) }}}{(1+|\Omega |t)^{\gamma _{1}\left( 1-\frac{2}{p}\right) }}%
\frac{1}{(1+|\Omega |t)^{\gamma _{2}\left( 1-\frac{2}{p}\right) }} \\
\leq (|\Omega |t)^{-\gamma _{2}\left( 1-\frac{2}{p}\right) },
\end{equation*}%
where $\gamma _{1},\gamma _{2}>0$, $\gamma _{1}+\gamma _{2}=\frac{1}{2}$ and
$\gamma _{2}<\frac{1}{2}\left( 1-\frac{1}{e}\right) $, we can estimate
(similarly to Lemma \ref{nonlinear1})
\begin{equation*}
\left\Vert \mathfrak{B}(u,u)-\mathfrak{B}(v,v)\right\Vert _{L^{\theta
}(0,\infty ;\dot{B}_{p,q}^{s+\epsilon })}\leq C|\Omega |^{-\alpha
_{0}-\gamma _{2}\left( 1-\frac{2}{p}\right) }(\Vert u\Vert _{L^{\theta
}(0,\infty ;\dot{B}_{p,q}^{s})}+\Vert v\Vert _{L^{\theta }(0,\infty ;\dot{B}%
_{p,q}^{s})})\Vert u-v\Vert _{L^{\theta }(0,\infty ;\dot{B}_{p,q}^{s})}
\end{equation*}%
which implies
\begin{equation*}
|\Omega |^{\alpha }\left\Vert \mathfrak{B}(u,u)-\mathfrak{B}(v,v)\right\Vert
_{L^{\theta }(0,\infty ;\dot{B}_{p,q}^{s+\epsilon })}\leq C|\Omega |^{\alpha
-\alpha _{0}-\gamma _{2}\left( 1-\frac{2}{p}\right) -2\beta _{0}}.
\end{equation*}%
Thus, for $\alpha <\alpha _{0}+2\beta _{0}+\gamma _{2}\left( 1-\frac{2}{p}%
\right) $, we obtain the property (\ref{asymp-201}).
\end{remark}

In what follows, we address the asymptotic behavior of solutions in the
critical case ($s=1/2$).

\begin{theorem}
\label{propertiescritical_u_and_v} Let $2\leq q\leq \infty $ and let $u$ and
$v$ be mild solutions of (\ref{NSC}) with initial data $u_{0}$ and $v_{0}$
in $\mathcal{F}_{0}$, respectively. Then, for all $\alpha \geq 0$
\begin{equation}
\lim_{|\Omega |\rightarrow \infty }|\Omega |^{\alpha }\Vert u-v\Vert
_{L^{4}(0,\infty ;\dot{B}_{3,q}^{\frac{1}{2}})}=0\ \ \ \text{if and only if}%
\ \ \ \lim_{|\Omega |\rightarrow \infty }|\Omega |^{\alpha }\Vert T_{\Omega
}(t)(u_{0}-v_{0})\Vert _{L^{4}(0,\infty ;\dot{B}_{3,q}^{\frac{1}{2}})}=0.
\label{asym-crit-2}
\end{equation}%
Moreover, for each $t>0$, we have that%
\begin{equation}
\lim_{|\Omega |\rightarrow \infty }|\Omega |^{\alpha }\Vert u(t)-v(t)\Vert _{%
\dot{B}_{2,q}^{\frac{1}{2}}}=0  \label{asym-crit-1}
\end{equation}%
provided that
\begin{equation}
\lim_{|\Omega |\rightarrow \infty }|\Omega |^{\alpha }\left( \Vert T_{\Omega
}(t)(u_{0}-v_{0})\Vert _{\dot{B}_{2,q}^{\frac{1}{2}}}+\Vert T_{\Omega
}(t)(u_{0}-v_{0})\Vert _{L^{4}(0,\infty ;\dot{B}_{3,q}^{\frac{1}{2}%
})}\right) =0.  \label{asymp-hip-1}
\end{equation}
\end{theorem}

\textbf{Proof. } By the proof of Theorem \ref{theorem3}, we know that $u\in
L^{4}(0,\infty ;\dot{B}_{3,q}^{\frac{1}{2}})$ with

\begin{equation*}
\Vert u\Vert _{L^{4}(0,\infty ;\dot{B}_{3,q}^{\frac{1}{2}})}\leq
2C_{1}\delta ,\text{ for all }|\Omega |\geq \tilde{\Omega},
\end{equation*}%
and similarly for $v$. Thus,
\begin{equation}
\sup_{|\Omega |\geq \tilde{\Omega}}\Vert u\Vert _{L^{4}(0,\infty ;\dot{B}%
_{3,q}^{\frac{1}{2}})}<2C_{1}\delta \text{ and }\sup_{|\Omega |\geq \tilde{%
\Omega}}\Vert v\Vert _{L^{4}(0,\infty ;\dot{B}_{3,q}^{\frac{1}{2}%
})}<2C_{1}\delta .  \label{aux-assint-10}
\end{equation}%
Next, we estimate
\begin{equation*}
\Vert u-v\Vert _{L^{4}(0,\infty ;\dot{B}_{3,q}^{\frac{1}{2}})}\leq \Vert
T_{\Omega }(t)(u_{0}-v_{0})\Vert _{L^{4}(0,\infty ;\dot{B}_{3,q}^{\frac{1}{2}%
})}+C_{2}(\Vert u\Vert _{L^{4}(0,\infty ;\dot{B}_{3,q}^{\frac{1}{2}})}+\Vert
v\Vert _{L^{4}(0,\infty ;\dot{B}_{3,q}^{\frac{1}{2}})})\Vert u-v\Vert
_{L^{4}(0,\infty ;\dot{B}_{3,q}^{\frac{1}{2}})}
\end{equation*}%
which yields
\begin{equation}
(1-4C_{1}C_{2}\delta )|\Omega |^{\alpha }\Vert u-v\Vert _{L^{4}(0,\infty ;%
\dot{B}_{3,q}^{\frac{1}{2}})}\leq |\Omega |^{\alpha }\Vert T_{\Omega
}(t)(u_{0}-v_{0})\Vert _{L^{4}(0,\infty ;\dot{B}_{3,q}^{\frac{1}{2}})},
\label{aux-assint-11}
\end{equation}%
where $C_{1},C_{2}$ and $\delta $ are as in the proof of Theorem \ref%
{theorem3}. Since $1-4C_{1}C_{2}\delta >0$ and the term on the right side
converges to zero (by hypothesis), it follows the \textquotedblleft \textit{%
if}\textquotedblright\ part in (\ref{asym-crit-2}). For the reverse, we
write (\ref{aux-1}) as
\begin{equation}
T_{\Omega }(t)(u_{0}-v_{0})=u-v-[\mathfrak{B}(u,u)(t)-\mathfrak{B}(v,v)(t)]
\label{aux-1-2}
\end{equation}%
\ and proceed similarly.

Next, we turn to (\ref{asym-crit-1}). Applying the $\dot{B}_{2,q}^{\frac{1}{2%
}}$-norm and using Lemma \ref{nonlinearestimativecritical}, we obtain
\begin{equation}
\Vert u(t)-v(t)\Vert _{\dot{B}_{2,q}^{\frac{1}{2}}}\leq \Vert T_{\Omega
}(t)(u_{0}-v_{0})\Vert _{\dot{B}_{2,q}^{\frac{1}{2}}}+C(\Vert u\Vert
_{L^{4}(0,\infty ;\dot{B}_{3,q}^{\frac{1}{2}})}+\Vert v\Vert
_{L^{4}(0,\infty ;\dot{B}_{3,q}^{\frac{1}{2}})})\Vert u-v\Vert
_{L^{4}(0,\infty ;\dot{B}_{3,q}^{\frac{1}{2}})},  \label{aux-assint-30}
\end{equation}%
for each $t>0$. Multiplying (\ref{aux-assint-30}) by $\left\vert \Omega
\right\vert ^{\alpha }$, letting $\left\vert \Omega \right\vert \rightarrow
\infty $, and using (\ref{aux-assint-10}), (\ref{asymp-hip-1}) and (\ref%
{asym-crit-2}), we get (\ref{asym-crit-1}). \fin

\begin{remark}
Notice that we can take $v_{0}=0$ and $v=0$ in Theorems \ref%
{properties_of_u_and_v_Omega} and \ref{propertiescritical_u_and_v} and
obtain asymptotic behavior properties for $u=u_{\Omega }$ as $\left\vert
\Omega \right\vert \rightarrow \infty $. In particular, in Theorem \ref%
{propertiescritical_u_and_v}, we have that
\begin{equation}
\lim_{|\Omega |\rightarrow \infty }|\Omega |^{\alpha }\Vert u_{\Omega }\Vert
_{L^{4}(0,\infty ;\dot{B}_{3,q}^{\frac{1}{2}})}=0\text{ provided that }%
\lim_{|\Omega |\rightarrow \infty }|\Omega |^{\alpha }\Vert T_{\Omega
}(t)u_{0}\Vert _{L^{4}(0,\infty ;\dot{B}_{3,q}^{\frac{1}{2}})}=0.
\label{aux-asy-50}
\end{equation}%
In the case $\alpha =0$, notice that the latter limit holds true for $%
u_{0}\in \dot{B}_{2,q}^{1/2}(\mathbb{R}^{3})$ with $1\leq q<4$ (see Lemma %
\ref{critical_linear_semigroup}) and for all $u_{0}\in \mathcal{F}_{0}$.
\end{remark}

\end{document}